\documentclass[a4paper,11pt]{article}
\usepackage{amssymb}
\usepackage{amsfonts}
\usepackage{amsmath}

\input{tcilatex}

\begin{document}

\title{Schubert presentation of the cohomology ring of flag manifolds $G/T$}
\author{Haibao Duan\thanks{%
Supported by 973 Program 2011CB302400 and NSFC 11131008.} \\
Institute of Mathematics, Chinese Academy of Sciences,\\
Beijing 100190, P. R. China\\
dhb@math.ac.cn \and Xuezhi Zhao\thanks{%
Supported by NSF of China (10931005), and a project of Beijing Municipal
Education Commission (PHR201106118).} \\
Department of Mathematics, Capital Normal University,\\
Beijing 100048, P. R. China\\
zhaoxve@mail.cnu.edu.cn}
\date{}
\maketitle

\begin{abstract}
Let $G$ be a compact connected Lie group with a maximal torus $T$. In the
context of Schubert calculus we present the integral cohomology ring $%
H^{\ast }(G/T)$ by a minimal system of generators and relations.

\begin{description}
\item \textsl{2000 Mathematical Subject Classification:} 14M15, 57T15.

\item \textsl{Key words:} Lie groups, flag manifolds, cohomology.
\end{description}
\end{abstract}

\section{Introduction}

Let $G$ be a compact connected Lie group with Lie algebra $L(G)$, and
exponential map $\exp :L(G)\rightarrow G$. For a non--zero vector $u\in L(G)$
the centralizer $P$ of the one parameter subgroup $\left\{ \exp (tu)\in
G\mid t\in \mathbb{R}\right\} $ on $G$ is a \textsl{parabolic subgroup} of $%
G $. The homogeneous space $G/P$ is canonically a projective variety, called
a \textsl{flag manifold} of $G$ \cite{Ku,LG}. If the vector $u$ is
non-singular the centralizer $P$ is a maximal torus $T$ on $G$, and the flag
manifold $G/T$ is also known as \textsl{the complete flag manifold} of $G$.

Schubert calculus \cite{Sch} begun with the intersection theory of the 19th
century, together with its applications to enumerative geometry. Clarifying
this calculus is an important problem of algebraic geometry \cite{CP,F,Kl,So}%
. van der Waerden and A. Weil, who secured the foundation of modern
intersection theory \cite{Wa}, attributed this calculus to the determination
of the integral cohomology ring $H^{\ast }(G/P)$ of flag manifolds $G/P$ 
\cite{Wa0}, \cite[p.331]{W}.

The cohomology of flag manifolds has now been well understood. The basis
theorem of Chevalley \cite{Ch, BGG} assures that the Schubert classes on a
flag manifold $G/P$ form an additive basis of the cohomology $H^{\ast }(G/P)$%
; an explicit formula for multiplying the basis elements was available by
the authors \cite{D,DZ1,DZ5}. However, concerning many relevant studies \cite%
{D0,D1,DZ3,DZ4,F,Ku,LG} such a description of the ring $H^{\ast }(G/P)$ is
not a practical one, as the number of Schubert classes on a flag manifold is
normally very large, not to mention the number of \textsl{the structure
constants} required to expand the products of Schubert classes. It is
natural to ask for a concise presentation of the ring $H^{\ast }(G/P)$ as
that characterized by the following notion.

Given a set $\left\{ x_{{1}},\cdots ,x_{{k}}\right\} $ of $k$ elements let $%
\mathbb{Z}[x_{{1}},\cdots ,x_{{k}}]$ be the ring of polynomials in $x_{{1}%
},\cdots ,x_{{k}}$ with integers coefficients. For a subset $\left\{
f_{1},\cdots ,f_{m}\right\} $ $\subset \mathbb{Z}[x_{{1}},\cdots ,x_{{k}}]$
write $\langle f_{1},\cdots ,f_{m}\rangle $ for the ideal generated by $%
f_{1},\cdots ,f_{m}$.

\bigskip

\noindent \textbf{Definition 1.1.} A \textsl{Schubert presentation} of the
integral cohomology ring of a flag manifold $G/P$ is an isomorphism

\begin{enumerate}
\item[(1.1)] $H^{\ast }(G/P)=\mathbb{Z}[x_{{1}},\cdots ,x_{{k}}]/\langle
f_{1},\cdots ,f_{m}\rangle $,
\end{enumerate}

\noindent where $\left\{ x_{{1}},\cdots ,x_{{k}}\right\} $ is a set of
Schubert classes on $G/P$ that generates the ring $H^{\ast }(G/P)$
multiplicatively, and where both the numbers $k$ and $m$ in (1.1) are
minimal.$\square $

\bigskip

Priory to the use of Schubert classes as generators, the numbers $k$ and $m$
in Definition 1.1 can be seen to be invariants of the ring $H^{\ast }(G/P)$.
Indeed, let $D(H^{\ast }(G/P))\subset H^{\ast }(G/P)$ be the ideal of
decomposable elements, and denote by $h(G,P)$ the cardinality of a basis for
the quotient group $H^{\ast }(G/P)/D(H^{\ast }(G/P))$. Then $k=h(G,P)-1$.
Furthermore, if one changes the generating set in (1.1) to $x_{{1}}^{\prime
},\cdots ,x_{{k}}^{\prime }$, then each old generator $x_{{i}}$ can be
expressed as a polynomial $g_{i}$ in the new ones $x_{{1}}^{\prime },\cdots
,x_{{k}}^{\prime }$. The invariance of the number $m$ is shown by the
presentation

\begin{quote}
$H^{\ast }(G/P)=\mathbb{Z}[x_{{1}}^{\prime },\cdots ,x_{{k}}^{\prime
}]/\langle f_{1}^{\prime },\cdots ,f_{m}^{\prime }\rangle $,
\end{quote}

\noindent where the polynomial $f_{j}^{\prime }$ is obtained from $f_{j}$ by
substituting the polynomial $g_{i}$ for $x_{{i}}$, $1\leq j\leq m$.

Among all the flag manifolds $G/P$ of a Lie group $G$ it is the complete
flag manifold $G/T$ that is of crucial importance. The inclusion $T\subset
P\subset G$ of subgroups induces the fibration $P/T\hookrightarrow G/T%
\overset{\pi }{\longrightarrow }G/P$ in which the induced map $\pi ^{\ast }$
embeds $H^{\ast }(G/P)$ as a subring of $H^{\ast }(G/T)$ (see Lemma 2.3). In
this paper we establish a Schubert presentation for the cohomologies of
complete flag manifolds $G/T$.

Recall that all the $1$--connected simple Lie groups consist of the three
infinite families $SU(n+1),Sp(n),$ $Spin(n+2)$, $n\geq 2,$ of classical
groups, as well as the five exceptional ones: $G_{2},F_{4},E_{6},E_{7},E_{8}$%
. It is also known that, for any compact connected Lie group $G$ with a
maximal torus $T$, one has a diffeomorphism $G/T=G_{1}/T_{1}\times \cdots
\times G_{k}/T_{k}$ \noindent with each $G_{i}$ a $1$--connected simple Lie
group and with $T_{i}\subset G_{i}$ a maximal torus. Moreover, by the basis
theorem of Chevalley (see Theorem 2.1) the integral cohomology $H^{\ast
}(G/T)$ is torsion free. Therefore, the problem of finding a Schubert
presentation of the ring $H^{\ast }(G/T)$ is reduced by the K\"{u}nneth
formula to the special cases where the group $G$ is $1$--connected and
simple. For this reason we assume in the remaining part of the paper that
the Lie group $G$ under consideration is simple. In addition, the
cohomologies are over the ring $\mathbb{Z}$ of integers, unless otherwise
stated.

For a Lie group $G$ of rank $n$ and let $\{\omega _{1},\cdots ,\omega
_{n}\}\subset H^{2}(G/T)$ be a set of fundamental dominant weights of $G$ 
\cite{BH}, and set $m=h(G,T)-n-1$. Our main result is

\bigskip

\noindent \textbf{Theorem 1.2.} \textsl{There exist a set }$\left\{
y_{d_{1}},\cdots ,y_{d_{m}}\right\} $\textsl{\ of }$m$ \textsl{Schubert
classes on }$G/T$\textsl{,\ with }$1<d_{1}<\cdots <d_{m}$ \textsl{and} $\deg
y_{d_{j}}=2d_{j}$\textsl{,} \textsl{so that the inclusion }$\left\{ \omega
_{1},\cdots ,\omega _{n},y_{d_{1}},\cdots ,y_{d_{m}}\right\} \subset $%
\textsl{\ }$H^{\ast }(G/T)$ \textsl{induces an isomorphism}

\begin{enumerate}
\item[(1.2)] $H^{\ast }(G/T)=\mathbb{Z}[\omega _{1},\cdots ,\omega
_{n},y_{d_{1}},\cdots ,y_{d_{m}}]/\left\langle
e_{i},f_{j},g_{j}\right\rangle _{1\leq i\leq k;1\leq j\leq m}$\textsl{,}
\end{enumerate}

\noindent \textsl{where}

\textsl{i)} $k=n-m$ \textsl{for all }$G\neq E_{8}$\textsl{\ but} $k=n-m+2$ 
\textsl{for} $G=E_{8}$\textsl{;}

\textsl{ii) }$e_{i}\in \left\langle \omega _{1},\cdots ,\omega
_{n}\right\rangle $\textsl{, }$1\leq i\leq k$\textsl{;}

\textsl{iii) the pair }$(f_{j},g_{j})$\textsl{\ of polynomials is related to
the class }$y_{d_{j}}$\textsl{\ in the fashion}

\begin{quote}
$f_{j}$\textsl{\ }$=$\textsl{\ }$p_{j}\cdot y_{d_{j}}+\alpha _{j}$\textsl{, }%
$g_{j}=y_{d_{j}}^{k_{j}}+\beta _{j}$\textsl{, }$1\leq j\leq m$\textsl{, }
\end{quote}

\noindent \textsl{with }$p_{j}\in \{2,3,5\}$\textsl{\ and }$\alpha
_{j},\beta _{j}\in \left\langle \omega _{1},\cdots ,\omega _{n}\right\rangle 
$\textsl{.}$\square $\noindent

\bigskip

A set $S=\left\{ y_{d_{1}},\cdots ,y_{d_{m}}\right\} $\textsl{\ }of Schubert
classes on\textsl{\ }$G/T$ satisfying (1.2) will be called \textsl{a set of
special Schubert classes} on\textsl{\ }$G/T$. In the course of showing
Theorem 1.2 a set of special Schubert classes, as well as the corresponding
system $\left\{ e_{i},f_{j},g_{j}\right\} $ of polynomials, will be made
explicit for each simple Lie group. Along the way an algebraic criterion for
a set of Schubert classes on\textsl{\ }$G/T$ to be special is given by
Theorem 6.3.

Since the set $\left\{ \omega _{1},\cdots ,\omega _{n}\right\} $ of
fundamental weights is precisely the Schubert basis on the group $H^{2}(G/T)$
\cite{DZ2}, the presentation (1.2) describes the ring $H^{\ast }(G/T)$ by
certain Schubert classes on $G/T$. It is worthwhile to know whether it is
indeed a Schubert presentation of the ring $H^{\ast }(G/T)$.

\bigskip

\noindent \textbf{Theorem 1.3. }\textsl{If }$G\neq E_{8}$\textsl{\ the
formula (1.2) is a Schubert presentation of }$H^{\ast }(G/T)$\textsl{. If }$%
G=E_{8}$ \textsl{a Schubert presentation of the ring }$H^{\ast }(E_{8}/T)$%
\textsl{\ is}

\begin{enumerate}
\item[(1.3)] $\mathbb{Z}[\omega _{1},\cdots ,\omega _{8},y_{d_{1}},\cdots
,y_{d_{7}}]/\left\langle e_{i},f_{j},g_{t},\phi \right\rangle _{1\leq i\leq
3;1\leq j\leq 7,t=1,2,3,5}$, \textsl{where}

\textsl{a) the Schubert classes }$y_{d_{1}},\cdots ,y_{d_{7}}$\textsl{\ and
the polynomials }$e_{i},f_{j},g_{t}$\textsl{\ are the same as that in (1.2)
for the case of }$G=E_{8}$\textsl{;}

\textsl{b) }$\phi =2y_{6}^{5}-y_{10}^{3}+y_{15}^{2}+\beta $\textsl{\ with }$%
\beta \in \left\langle \omega _{1},\cdots ,\omega _{8}\right\rangle $\textsl{%
.}$\square $
\end{enumerate}

This paper is arranged as follows. Section 2 develops cohomology properties
for firations in flag manifolds. Granted with the packages "\textsl{the Chow
ring of Grassmannians}" and "\textsl{Giambelli polynomials}" compiled in%
\textsl{\ }\cite[Section 2.6]{DZ2} initial data facilitating our computation
are generated in Sections 3 and 4. With these preparations the presentation
(1.2) for the exceptional Lie groups are obtained in Section 5. Finally,
Theorems 1.2 and 1.3 are established in Section 6.

Certain relations on the ring $H^{\ast }(G/T)$ may be seen as detailed.
However, they are useful for encoding the topology of the corresponding Lie
group $G$. Using the set $\left\{ e_{i},f_{j},g_{j}\right\} $ of polynomials
in (1.2) one can construct uniformly the integral cohomology of compact Lie
groups \cite{D0,DZ3}, deduce explicit formulae for the generalized Weyl
invariants of $G$ in a characteristic $p$ \cite[Propositions 5.5--5.7]{DZ4},
and determine the structure of the $\func{mod}p$ cohomology $H^{^{\ast }}(G;%
\mathbb{F}_{p})$ as a module over the Steenrod algebra $\mathcal{A}_{p}$ 
\cite{DZ4}.

\section{Fibrations in flag manifolds}

Let $G$ be Lie group with maximal torus $T$ and Cartan subalgebra $L(T)$.
Equip the Lie algebra $L(G)$ with an inner product $(,)$ so that the adjoint
representation acts as isometries on $L(G)$. Assume that the rank of $G$ is $%
n=\dim T$, and a system $\{\beta _{1},\cdots ,\beta _{n}\}$ of simple roots
of $G$ is so ordered as the vertices in the Dynkin diagram of $G$ pictured
in \cite[p.58]{H}. Then the \textsl{Weyl group} of $G$ is the subgroup $%
W\subset Aut(L(T))$ generated by the reflections $\sigma _{i}$ in the
hyperplanes $L_{i}\subset L(T)$ perpendicular to the roots $\beta _{i}$, $%
1\leq i\leq n$. By the relation $H^{2}(G/T)\otimes \mathbb{R}=L(T)$ due to
Borel and Hirzebruch \cite{BH} the set $\{\omega _{{1}},\ldots ,\omega _{{n}%
}\}\subset H^{2}(G/T)$ of \textsl{fundamental dominant weights} of $G$ can
be regarded as the basis of the space $L(T)$ defined by the formulae

\begin{quote}
$2(\beta _{{i}},\omega _{{j}})/(\beta _{{i}},\beta _{{i}})=\delta _{{i,j}}$, 
$1\leq i,j\leq n$.
\end{quote}

For a parabolic subgroup $P$ on $G$ let $W$ and $W^{\prime }$ be the Weyl
groups of $G$ and $P$, respectively. In term of the length function $l:$ $%
W\rightarrow \mathbb{Z}$ on $W$ the set $\overline{W}$ of left cosets of $%
W^{\prime }$ in $W$ can be identified with the subset of $W$

\begin{quote}
$\overline{W}=\{w\in W\mid l(w_{{1}})\geq l(w)$, $w_{{1}}\in wW^{\prime }\}$,
\end{quote}

\noindent see \cite[5.1]{BGG}. It follows that every element $w\in \overline{%
W}$ admits a decomposition $w=\sigma _{{i}_{{1}}}\circ \cdots \circ \sigma _{%
{i}_{{r}}}$ with $1\leq i_{{1}},\ldots ,i_{{r}}\leq n$ and $r=$ $l(w)$. This
decomposition is called \textsl{minimized}, written $w=\sigma \lbrack i_{{1}%
},\ldots ,i_{{r}}]$, if the relation $(i_{{1}},\ldots ,i_{{r}})\leq (j_{{1}%
},\ldots ,j_{{r}})$ holds for any $(j_{{1}},\ldots ,j_{{r}})$ satisfying $%
w=\sigma _{{j}_{{1}}}\circ \cdots \circ \sigma _{{j}_{{r}}}$, where $\leq $
means the lexicographical order on multi--indices.

For an element $w\in \overline{W}$ with minimized decomposition $\sigma
\lbrack i_{{1}},\ldots ,i_{{r}}]$ the \textsl{Schubert variety} $X_{w}$
associated to $w$ is the image of the composition

\begin{quote}
$K_{{i}_{{1}}}\times \cdots \times K_{{i}_{{r}}}\rightarrow G\overset{p}{%
\rightarrow }G/P$, $(k_{{1}},\ldots ,k_{{r}})\longmapsto p(k_{{1}}\cdot
\cdots \cdot k_{{r}})$,
\end{quote}

\noindent where $K_{{i}}\subset G$ is the centralizer of $\exp (L_{{i}})$ in 
$G$, $p$ is the obvious quotient map, and where the product $\cdot $ takes
place in $G$. In \cite{Ch} Chevalley announced the following remarkable
cellular decomposition on the flag manifold $G/P$

\begin{enumerate}
\item[(2.1)] $G/P=\underset{w\in \overline{W}}{\bigcup }X_{w}$, $\quad \dim
_{\mathbb{R}}X_{w}=2l(w)$ (see also \cite{BGG,De}).
\end{enumerate}

\noindent The \textsl{Schubert class }$s_{w}\in H^{2l(w)}(G/P)$
corresponding to $w\in \overline{W}$ is defined to be the Kronecker dual of
the fundamental classes\textsl{\ }$[X_{w}]\in H_{2l(w)}(G/P)$. Since only
even dimensional cells are involved in the decomposition (2.1) one has the
next result, called the \textsl{basis theorem} of Schubert calculus.

\bigskip

\noindent \textbf{Theorem 2.1 (see \cite{Ch,BGG,De})}. \textsl{The set of
Schubert classes }$\{s_{w}\in H^{\ast }(G/P)\mid $\textsl{\ }$w\in \overline{%
W}\}$\textsl{\ constitutes a basis of the graded group }$H^{\ast }(G/P)$%
\textsl{.}$\square $

\bigskip

For a subset $I\subseteq \{1,\cdots ,n\}$ let $P_{I}$ be the centralizer of
the $1$--parameter subgroup $\alpha :\mathbb{R}\rightarrow G$, $\alpha
(t)=\exp (t\sum\limits_{i\in I}\omega _{i})$ on $G$. Useful information on
the geometry of the flag manifold $G/P_{I}$ is

\bigskip

\noindent \textbf{Lemma 2.2}(\cite{DZ2,DZ5})\textbf{. }\textsl{The
centralizer of any }$1$\textsl{--parameter subgroup on }$G$\textsl{\ is
conjugate\ to a subgroup }$P_{I}$\textsl{\ for some }$I\subseteq \{1,\cdots
,n\}$\textsl{. Moreover,}

\textsl{i) }$P_{I}$\textsl{\ is a parabolic subgroup\ its Dynkin diagram is
obtained from that of }$G$\textsl{\ by deleting the vertices }$\beta _{i}$%
\textsl{\ with }$i\in I$\textsl{,} \textsl{and the edges adjoining to it;}

\textsl{ii) the Weyl group }$W_{I}$\textsl{\ of }$P_{I}$\textsl{\ is the
subgroup of }$W$\textsl{\ generated by }$\sigma _{{j}},j\notin I$\textsl{;}

\textsl{iii)} \textsl{the Schubert basis of }$H^{\ast }(G/P_{I})$ \textsl{is 
}$\{s_{w}\mid w\in W/W_{I}\}.\square $

\bigskip

Property i) of Lemma 2.2 characterizes the subgroup $P_{I}$ only up to its
local type. A method to decide the isomorphism type of $P_{I}$ is given in 
\cite{DL}.

For a proper subset $I\subset \{1,\cdots ,n\}$ the inclusion $T\subset P_{I}$
$\subset G$ of subgroups induces the fibration in flag manifolds

\begin{enumerate}
\item[(2.2)] $P_{I}/T\overset{i}{\hookrightarrow }G/T\overset{\pi }{%
\rightarrow }G/P_{I}$.
\end{enumerate}

\noindent The next result implies that the cohomologies of the fiber space $%
P_{I}/T$ and the base space $G/P_{I}$ are much simpler than that of the
total space $G/T$.

\bigskip

\noindent \textbf{Lemma 2.3.} \textsl{With respect to the inclusion }$%
W_{I}\subset W$\textsl{\ the induced map }$i^{\ast }$\textsl{\ identifies
the subset }$\{s_{w}\}_{w\in W_{I}\subset W}$\textsl{\ of the Schubert basis
of }$H^{\ast }(G/T)$\textsl{\ with the Schubert basis }$\{s_{w}\}_{w\in
W_{I}}$\textsl{\ of }$H^{\ast }(P_{I}/T)$.

\textsl{With respect to the inclusion }$W/W_{I}\subset W$\textsl{\ the
induced map }$\pi ^{\ast }$ \textsl{identifies the Schubert basis }$%
\{s_{w}\}_{w\in W/W_{I}}$\textsl{\ of }$H^{\ast }(G/P_{I})$\textsl{\ with
the subset }$\{s_{w}\}_{w\in W/W_{I}}$\textsl{\ of the Schubert basis }$%
\{s_{w}\}_{w\in W}$ \textsl{of }$H^{\ast }(G/T)$\textsl{.}

\bigskip

\noindent \textbf{Proof.} These come directly from the next two properties
of Schubert varieties (e.g. \cite[Section 2]{DZ2}). With respect to the cell
decompositions (2.1) on the three flag manifolds $P_{I}/T$, $G/T$ and $%
G/P_{I}$ one has:

i) for each $w\in W_{I}\subset W$ the fiber inclusion $i$ carries the
Schubert variety $X_{w}$ on $P_{I}/T$ identically onto the Schubert variety $%
X_{w}$ on $G/T$;

ii) for each $w\in W/W_{I}\subset W$ the projection $\pi $ restricts to a
degree $1$ map from the Schubert variety $X_{w}$ on $G/T$ to the
corresponding Schubert variety on $G/P_{I}$.$\square $

\bigskip

\noindent \textbf{Convention 2.4.} In view of Lemma 2.3 and for the
notational convenience, we shall make no difference in notation between an
element in $H^{\ast }(G/P_{I})$\ and its $\pi ^{\ast }$\ image in $H^{\ast
}(G/T)$, and between a Schubert class on $P_{I}/T$\ and its $i^{\ast }$\
pre-image on $G/T$.$\square $

\bigskip

To formulate the ring $H^{\ast }(G/T)$ in question from the simpler ones $%
H^{\ast }(P_{I}/T)$ and $H^{\ast }(G/P_{I})$ assume that $\{y_{1},\cdots
,y_{n_{1}}\}$\ is a subset of Schubert classes on $P_{I}/T$, $\{x_{1},\cdots
,x_{n_{2}}\}$\ is a subset of Schubert classes on $G/P_{I}$, and that with
respect to them one has the following presentations of the cohomologies

\begin{enumerate}
\item[(2.3)] $H^{\ast }(P_{I}/T)=\frac{\mathbb{Z}[y_{i}]_{1\leq i\leq n_{1}}%
}{\left\langle h_{s}\right\rangle _{1\leq s\leq m_{1}}}$; $H^{\ast
}(G/P_{I})=\frac{\mathbb{Z}[x_{j}]_{1\leq j\leq n_{2}}}{\left\langle
r_{t}\right\rangle _{1\leq t\leq m_{2}}}$,
\end{enumerate}

\noindent where $h_{s}\in \mathbb{Z}[y_{i}]_{1\leq i\leq n_{1}}$, $r_{t}\in 
\mathbb{Z}[x_{j}]_{1\leq j\leq n_{2}}$.

\bigskip

\noindent \textbf{Lemma 2.5.} \textsl{The inclusions }$y_{i},x_{j}\in
H^{\ast }(G/T)$ \textsl{induces a surjective map}

\begin{quote}
\textsl{\ }$\varphi :\mathbb{Z}[y_{i},x_{j}]_{1\leq i\leq n_{1},1\leq j\leq
n_{2}}\rightarrow H^{\ast }(G/T)$\textsl{. }
\end{quote}

\noindent \textsl{Furthermore, if }$\{\rho _{s}\}_{1\leq s\leq m_{1}}\subset 
\mathbb{Z}[y_{i},x_{j}]$\textsl{\ is a system satisfying}

\begin{enumerate}
\item[(2.4)] \textsl{\ }$\varphi (\rho _{s})=0$\textsl{\ and }$\rho _{s}\mid
_{x_{j}=0}=h_{s}$\textsl{,}
\end{enumerate}

\noindent \textsl{then }$\varphi $\textsl{\ induces a ring isomorphism}

\begin{enumerate}
\item[(2.5)] $H^{\ast }(G/T)=\mathbb{Z}[y_{i},x_{i}]_{1\leq i\leq
n_{1},1\leq j\leq n_{2}}/\left\langle \rho _{s},r_{t}\right\rangle _{1\leq
s\leq m_{1},1\leq t\leq m_{2}}$\textsl{.}
\end{enumerate}

\noindent \textbf{Proof.} Lemma 2.3, together with Convention 2.4, implies
that the map $\varphi $ is surjective. It remains to show that for a $g\in 
\mathbb{Z}[y_{i},x_{j}]_{1\leq i\leq n_{1},1\leq j\leq n_{2}}$ the relation $%
\varphi (g)=0$ implies that $g\in \left\langle \rho _{s},r_{t}\right\rangle
_{1\leq s\leq m_{1},1\leq t\leq m_{2}}$. By Lemma 2.3 and by the
Leray--Hirsch property \cite[p.231]{Hus} of the fibration (2.2) one has the
following presentation of $H^{\ast }(G/T)$ a module over its subring $%
H^{\ast }(G/P_{I})$

\begin{quote}
$H^{\ast }(G/T)=H^{\ast }(G/P_{I})\{1,s_{w}\}_{w\in W_{I}}$.
\end{quote}

\noindent It follows from the presentation of the ring $H^{\ast }(P_{I}/T)$
in (2.3) and the assumption (2.4) that, for any polynomial $g\in \mathbb{Z}%
[y_{i},x_{j}]$ one has

\begin{quote}
$g\equiv \sum\nolimits_{w\in W_{I}}g_{w}\cdot s_{w}$ mod $\left\langle \rho
_{s}\right\rangle _{1\leq s\leq m_{1}}$ with $g_{w}\in \mathbb{Z}%
[x_{j}]_{1\leq j\leq n_{2}}$.
\end{quote}

\noindent From this we find that $\varphi (g)=0$ implies that $\varphi
(g_{w})=0$, $w\in W_{I}$. That is $g_{w}\in \left\langle r_{t}\right\rangle
_{1\leq t\leq m_{2}}$, $w\in W_{I}$, by the presentation of the ring $%
H^{\ast }(G/P_{I})$ in (2.3). This completes the proof.$\square $

\section{Cohomology of generalized Grassmannians}

If $I=\{k\}$ is a singleton the flag manifold $G/P_{\{k\}}$ is called the 
\textsl{Grassmannians of }$G$\textsl{\ corresponding to the weight }$\omega
_{k}$ \cite{DZ2}. Using the table below we associate to each exceptional Lie
group $G$ a Grassmannian $G/P_{\{k\}}$

\begin{enumerate}
\item[(3.1)] 
\begin{tabular}{l|lllll}
\hline\hline
$G$ & $G_{2}$ & $F_{4}$ & $E_{6}$ & $E_{7}$ & $E_{8}$ \\ \hline
$k$ & $1$ & $1$ & $2$ & $2$ & $2$ \\ \hline
$P_{\{k\}}$ & $SU(2)\cdot S^{1}$ & $Sp(3)\cdot S^{1}$ & $SU(6)\cdot S^{1}$ & 
$SU(7)\cdot S^{1}$ & $SU(8)\cdot S^{1}$ \\ \hline
$P_{\{k\}}^{s}$ & $SU(2)$ & $Sp(3)$ & $SU(6)$ & $SU(7)$ & $SU(8)$ \\ 
\hline\hline
\end{tabular}%
,
\end{enumerate}

\noindent where in the third row the subgroups $P_{\{k\}}$ are presented by
their local types determined by i) of Lemma 2.2, and where the group $%
P_{\{k\}}^{s}$ in the fourth row is the simple part of the group $P_{\{k\}}$%
. Our approach to the ring $H^{\ast }(G/T)$ amounts to apply Lemma 2.5 to
the fibration in flag manifolds

\begin{enumerate}
\item[(3.2)] $P_{\{k\}}/T\overset{i}{\hookrightarrow }G/T\overset{\pi }{%
\rightarrow }G/P_{\{k\}}$.
\end{enumerate}

\noindent To this end an account for the cohomologies of the base spaces $%
G/P_{\{k\}}$ is required.

In \cite{DZ2} a program calculating Schubert presentation of a Grassmannian $%
G/P_{\{k\}}$ has been compiled, whose function is briefed below:

\begin{quote}
\textbf{Algorithm:} \textsl{The Chow ring of Grassmannians}

\textbf{Input:} \textsl{The Cartan matrix }$C=(c_{ij})_{n}$\textsl{\ of }$G$%
\textsl{, }and an integer $k\in \left\{ 1,...,n\right\} $;

\textbf{Output:} \textsl{A Schubert presentation of the cohomology} $H^{\ast
}(G/P_{\{k\}})$.
\end{quote}

\noindent As applications of the algorithm, Schubert presentations for the
five Grassmannians tabulated in (3.1) have been obtained. To state the
results we introduce for each of the Grassmannians $G/P_{\{k\}}$ a set of
Schubert classes $\left\{ s_{w}\right\} $ on $G/P_{\{k\}}$ in terms of the
minimized decomposition $\sigma \lbrack i_{{1}},\ldots ,i_{{r}}]$ of the
corresponding $w$, together with their abbreviations $y_{i}$'s (with $\deg
y_{i}=2i$), in the table below:

\begin{enumerate}
\item[(3.3)] 
\begin{tabular}{|l|l|l|l|}
\hline
$y_{i}$ & $G_{2}/P_{\{1\}}$ & $F_{4}/P_{\{1\}}$ & $E_{n}/P_{\{2\}},\text{ }%
n=6,7,8$ \\ \hline\hline
$y_{3}$ & $\sigma _{\lbrack 1,2,1]}$ & $\sigma _{\lbrack 3,2,1]}$ & $\sigma
_{\lbrack 5,4,2]}\text{, }n=6,7,8$ \\ 
$y_{4}$ &  & $\sigma _{\lbrack 4,3,2,1]}$ & $\sigma _{\lbrack 6,5,4,2]}\text{%
, }n=6,7,8$ \\ 
$y_{5}$ &  &  & $\sigma _{\lbrack {7,6,5,4,2}]}\text{, }n=7,8$ \\ 
$y_{6}$ &  & $\sigma _{\lbrack 3,2,4,3,2,1]}$ & $\sigma _{\lbrack
1,3,6,5,4,2]}\text{, }n=6,7,8$ \\ 
$y_{7}$ &  &  & $\sigma _{\lbrack 1,3,7,6,5,4,2]}$, $n=7,8$ \\ 
$y_{8}$ &  &  & $\sigma _{\lbrack 1,3,8,7,6,5,4,2]}$, $n=8$ \\ 
$y_{9}$ &  &  & $\sigma _{\lbrack 1,5,4,3,7,6,5,4,2]}$,$\text{ }n=7,8$ \\ 
$y_{10}$ &  &  & $\sigma _{\lbrack {1,6,5,4,3,7,6,5,4,2}]}\text{, }n=8$ \\ 
$y_{15}$ &  &  & $\sigma _{\lbrack 5,4,2,3,1,6,5,4,3,8,7,6,5,4,2]}\text{, }%
n=8$ \\ \hline
\end{tabular}%
.
\end{enumerate}

\noindent \textbf{Theorem 3.1. }\textsl{With respect to the special Schubert
classes on }$G/P_{\{k\}}$\textsl{\ given in (3.3), the Schubert
presentations of the integral cohomology rings of }$G/P_{\{k\}}$\textsl{\
are:}

\begin{enumerate}
\item[(3.4)] $H^{\ast }(G_{2}/SU(2)\cdot S^{1})=\mathbb{Z}[\omega
_{1},y_{3}]/\left\langle r_{3},r_{6}\right\rangle $\textsl{, where}

$r_{3}=2y_{3}-\omega _{1}^{3}$;$\quad $

$r_{6}=y_{3}^{2}$.

\item[(3.5)] $H^{\ast }(F_{4}/Sp(3)\cdot S^{1})=\mathbb{Z}[\omega
_{1},y_{3},y_{4},y_{6}]/\left\langle r_{3},r_{6},r_{8},r_{12}\right\rangle $%
\textsl{, where}

$r_{3}=2y_{3}-\omega _{1}^{3}$;$\qquad $

$r_{6}=2y_{6}+y_{3}^{2}-3\omega _{1}^{2}y_{4}$;

$r_{8}=3y_{4}^{2}-\omega _{1}^{2}y_{6}$;

$r_{12}=y_{6}^{2}-y_{4}^{3}.$

\item[(3.6)] $H^{\ast }(E_{6}/SU(6)\cdot S^{1})=\mathbb{Z}[\omega
_{2},y_{3},y_{4},y_{6}]/\left\langle r_{6},r_{8},r_{9},r_{12}\right\rangle $%
\textsl{, where}

$r_{6}=2y_{6}+y_{3}^{2}-3\omega _{2}^{2}y_{4}+2\omega _{2}^{3}y_{3}-\omega
_{2}^{6};$

$r_{8}=3y_{4}^{2}-6\omega _{2}y_{3}y_{4}+\omega _{2}^{2}y_{6}+5\omega
_{2}^{2}y_{3}^{2}-2\omega _{2}^{5}y_{3};$

$r_{9}=2y_{3}y_{6}-\omega _{2}^{3}y_{6};$

$r_{12}=y_{6}^{2}-y_{4}^{3}.$

\item[(3.7)] $H^{\ast }(E_{7}/SU(7)\cdot S^{1})=\mathbb{Z}[\omega
_{2},y_{3},y_{4},y_{5},y_{6},y_{7},y_{9}]/$\textsl{\ }$\left\langle
r_{j}\right\rangle _{j\in \mathcal{R}(7)}$ \textsl{with }$\mathcal{R}%
(7)=\{6,8,9,10,12,14,18\}$\textsl{,} \textsl{where}

$r_{6}=2y_{6}+y_{3}^{2}+2\omega _{2}y_{5}-3\omega _{2}^{2}y_{4}+2\omega
_{2}^{3}y_{3}-\omega _{2}^{6}$;

$r_{8}=$ \ $3y_{4}^{2}-2y_{3}y_{5}+2\omega _{2}y_{7}-6\omega
_{2}y_{3}y_{4}+\omega _{2}^{2}y_{6}+5\omega _{2}^{2}y_{3}^{2}+2\omega
_{2}^{3}y_{5}-2\omega _{2}^{5}y_{3}$;

$r_{9}=2y_{9}+2y_{4}y_{5}-2y_{3}y_{6}-4\omega _{2}y_{3}y_{5}\ -\omega
_{2}^{2}y_{7}+\omega _{2}^{3}y_{6}\ +2\omega _{2}^{4}y_{5}$;

$r_{10}=y_{5}^{2}-2y_{3}y_{7}+\omega _{2}^{3}y_{7}$;

$r_{12}=y_{6}^{2}+2y_{5}y_{7}-y_{4}^{3}+2y_{3}y_{9}+2y_{3}y_{4}y_{5}+2\omega
_{2}y_{5}y_{6}-6\omega _{2}y_{4}y_{7}+\omega _{2}^{2}y_{5}^{2}$;

$r_{14}=y_{7}^{2}-2y_{5}y_{9}+y_{4}y_{5}^{2}$;

$%
r_{18}=y_{9}^{2}+2y_{5}y_{6}y_{7}-y_{4}y_{7}^{2}-2y_{4}y_{5}y_{9}+2y_{3}y_{5}^{3}-\omega _{2}y_{5}^{2}y_{7} 
$.

\item[(3.8)] $H^{\ast }(E_{8}/SU(8)\cdot S^{1})=\mathbb{Z}[\omega
_{2},y_{3},y_{4},y_{5},y_{6},y_{7},y_{8},y_{9},y_{10},y_{15}]/\left\langle
r_{j}\right\rangle _{j\in \mathcal{R}(8)}$ \textsl{with }$\mathcal{R}%
(8)=\{8,9,10,12,14,15,18,20,24,30\}$\textsl{,} \textsl{where}

$r_{8}=3{y_{8}}-3{y_{4}^{2}}+2{y_{3}y_{5}}-2\omega _{2}{y_{7}}+6\omega _{2}{%
y_{3}y_{4}}-\omega _{2}^{2}{y_{6}-}\omega _{2}^{2}z_{6}-5\omega _{2}^{2}{%
y_{3}^{2}}$

$\qquad -2\omega _{2}^{3}{y_{5}}+2\omega _{2}^{5}{y_{3}}$;

$r_{9}=2y_{9}+2y_{4}y_{5}-2y_{3}y_{6}+\ \omega _{2}y_{8}-4\omega
_{2}y_{3}y_{5}\ -\omega _{2}^{2}y_{7}+\omega _{2}^{3}y_{6}\ +2\omega
_{2}^{4}y_{5}$;

$r_{10}=3{{y_{10}}}-2{y_{5}^{2}}+6{{y_{4}z_{6}}}-2{y_{3}y_{7}-}4\omega _{2}{%
y_{3}z}_{6}-\omega _{2}^{2}{y_{8}}+\omega _{2}^{3}{y_{7}}+2\omega
_{2}^{4}z_{6}$;

$r_{12}={y_{4}^{3}-y_{6}^{2}}+4{z}_{6}^{2}-2{y_{5}y_{7}}-6{y_{4}y_{8}}-2{%
y_{3}y_{9}}-2{y_{3}y_{4}y_{5}}-2{y_{3}^{2}z}_{6}-2\omega _{2}{y_{5}y_{6}}$

$\qquad -4\omega _{2}{y_{5}z}_{6}+6\omega _{2}{y_{4}y_{7}}+6\omega _{2}^{2}{%
y_{4}z}_{6}-2\omega _{2}^{2}{y_{3}y_{7}}+\omega _{2}^{5}{y_{7}}$;

$r_{14}={y_{7}^{2}}-6z_{6}{y_{8}}-3{y_{6}y_{8}}-2{y_{5}y_{9}}+3{y_{4}y{_{10}}%
+6y{_{4}^{2}z_{6}}}-2{y_{4}y_{5}^{2}}-4{y_{3}y_{5}}z_{6}$

$\qquad -4{y_{3}^{2}y_{8}}+4\omega _{2}{z_{6}y_{7}}+4\omega _{2}^{2}{z}%
_{6}^{2}+2\omega _{2}^{3}{y_{3}y_{8}}$;

$r_{15}=2y_{15}-y_{7}y_{8}+2y_{5}{y_{10}-2y_{5}^{3}+y_{4}y_{5}z_{6}+2}%
y_{3}z_{6}^{2}-2y_{3}y_{4}y_{8}+2\omega _{2}z_{6}y_{8}$

$\qquad -2\omega _{2}^{3}y_{6}^{2}+\omega _{2}^{3}y_{4}y_{8}$;

$r_{18}=y_{9}^{2}+9y_{10}y_{8}-3y_{10}y_{7}\omega _{2}-6y_{10}z_{6}\omega
_{2}^{2}-3y_{10}y_{5}\omega _{2}^{3}+4y_{9}z_{6}y_{3}$

$\qquad -2y_{9}y_{5}y_{4}-2y_{8}^{2}\omega
_{2}^{2}+y_{8}y_{7}y_{3}+18y_{8}z_{6}y_{4}+2y_{8}z_{6}y_{3}\omega
_{2}-2y_{8}y_{5}^{2}$

$\qquad -6y_{8}y_{5}y_{4}\omega _{2}+6y_{8}y_{5}y_{3}\omega
_{2}^{2}-y_{8}y_{5}\omega
_{2}^{5}-2y_{8}y_{4}y_{3}^{2}+y_{8}y_{4}y_{3}\omega _{2}^{3}-y_{7}^{2}y_{4}$

$\qquad +5y_{7}z_{6}y_{5}-10y_{7}z_{6}y_{4}\omega
_{2}+3y_{7}y_{5}y_{4}\omega
_{2}^{2}+3y_{7}y_{5}y_{3}^{2}-2y_{7}y_{5}y_{3}\omega _{2}^{3}$

$\qquad -10z_{6}^{3}+10z_{6}^{2}y_{5}\omega _{2}-16z_{6}^{2}y_{4}\omega
_{2}^{2}+4z_{6}^{2}y_{3}^{2}-4z_{6}y_{5}y_{4}y_{3}$

$\qquad +6z_{6}y_{5}y_{3}^{2}\omega _{2}-4z_{6}y_{5}y_{3}\omega
_{2}^{4}+2y_{5}^{3}\omega _{2}^{3}$;

$r_{20}=3(y_{10}+2y_{4}z_{6}-y_{5}^{2})^{2}+y_{8}(-4y_{9}y_{3}+2y_{9}\omega
_{2}^{3}-3y_{8}y_{4}+2y_{8}y_{3}\omega _{2}$

$\qquad -y_{8}\omega _{2}^{4}+y_{7}y_{5}+2z_{6}y_{5}\omega
_{2}-8z_{6}y_{3}^{2}+4z_{6}y_{3}\omega
_{2}^{3}+2y_{5}y_{4}y_{3}-y_{5}y_{4}\omega _{2}^{3})$;

$r_{24}=5z_{6}^{4}+y_{8}(6y_{10}z_{6}-12y_{10}y_{6}-18y_{10}y_{5}\omega
_{2}-12y_{10}y_{3}^{2}-5y_{9}y_{7}$

$\qquad -4y_{9}z_{6}\omega _{2}-2y_{9}y_{6}\omega _{2}+2y_{9}y_{5}\omega
_{2}^{2}-6y_{9}y_{4}y_{3}+3y_{9}y_{4}\omega _{2}^{3}-2y_{9}y_{3}^{2}\omega
_{2}$

$\qquad +4y_{8}^{2}-14y_{8}y_{7}\omega _{2}-6y_{8}z_{6}\omega
_{2}^{2}+6y_{8}y_{6}\omega _{2}^{2}-14y_{8}y_{5}y_{3}-21y_{8}y_{4}^{2}$

$\qquad +22y_{8}y_{4}y_{3}\omega _{2}-4y_{8}y_{4}\omega
_{2}^{4}+5y_{7}^{2}\omega _{2}^{2}+10y_{7}z_{6}\omega
_{2}^{3}-4y_{7}y_{5}y_{4}$

$\qquad +6y_{7}y_{5}y_{3}\omega _{2}+21y_{7}y_{4}^{2}\omega
_{2}-10y_{7}y_{4}y_{3}\omega _{2}^{2}+4y_{7}y_{4}\omega
_{2}^{5}+2y_{7}y_{3}^{3}$

$\qquad +4z_{6}^{2}y_{3}\omega
_{2}+15z_{6}y_{6}y_{4}-12z_{6}y_{6}y_{3}\omega
_{2}-12z_{6}y_{5}^{2}+4z_{6}y_{5}y_{3}\omega _{2}^{2}$

$\qquad -24z_{6}y_{4}^{2}\omega
_{2}^{2}-12z_{6}y_{4}y_{3}^{2}+30z_{6}y_{4}y_{3}\omega
_{2}^{3}-14z_{6}y_{4}\omega _{2}^{6}-2y_{6}^{2}y_{3}\omega _{2}$

$\qquad +4y_{6}y_{5}^{2}\ -6y_{6}y_{5}y_{3}\omega
_{2}^{2}+y_{6}y_{4}y_{3}^{2}+\ 12y_{5}^{3}\omega _{2}+2y_{4}^{4})$;

$%
r_{30}=z_{6}^{5}-(y_{10}+2z_{6}y_{4}-y_{5}^{2})^{3}+(y_{15}+y_{10}y_{5}+z_{6}^{2}y_{3}-z_{6}^{2}\omega _{2}^{3}+2z_{6}y_{5}y_{4}-y_{5}^{3})^{2} 
$

$\qquad +y_{8}(y_{15}y_{7}-4y_{15}z_{6}\omega
_{2}+4y_{15}y_{4}y_{3}-3y_{15}y_{4}\omega _{2}^{3}+6y_{10}y_{9}y_{3}$

$\qquad -3y_{10}y_{9}\omega
_{2}^{3}-9y_{10}y_{8}y_{4}+12y_{10}y_{8}y_{3}\omega _{2}-9y_{10}y_{8}\omega
_{2}^{4}-2y_{10}y_{7}y_{5}$

$\qquad -24y_{10}z_{6}^{2}+48y_{10}z_{6}y_{6}+50y_{10}z_{6}y_{5}\omega
_{2}-48y_{10}z_{6}y_{4}\omega _{2}^{2}+24y_{10}z_{6}y_{3}^{2}$

$\qquad -6y_{10}z_{6}y_{3}\omega
_{2}^{3}-2y_{10}y_{5}y_{4}y_{3}+6y_{9}y_{8}y_{5}-4y_{9}y_{8}y_{4}\omega
_{2}+8y_{9}y_{8}y_{3}\omega _{2}^{2}$

$\qquad -y_{9}y_{8}\omega
_{2}^{5}+5y_{9}y_{7}z_{6}-2y_{9}y_{7}y_{6}-2y_{9}y_{7}y_{5}\omega
_{2}+3y_{9}y_{7}y_{4}\omega _{2}^{2}$

$\qquad -5y_{9}y_{7}y_{3}^{2}+2y_{9}y_{7}y_{3}\omega _{2}^{3}\
-6y_{9}z_{6}^{2}\omega _{2}-4y_{9}z_{6}y_{6}\omega
_{2}+44y_{9}z_{6}y_{5}\omega _{2}^{2}$

$\qquad +4y_{9}z_{6}y_{4}y_{3}-2y_{9}z_{6}y_{3}^{2}\omega
_{2}+81y_{8}^{2}z_{6}+14y_{8}^{2}y_{6}\ -11y_{8}^{2}y_{5}\omega _{2}$

$\qquad +16y_{8}^{2}y_{4}\omega
_{2}^{2}+11y_{8}^{2}y_{3}^{2}-4y_{8}^{2}y_{3}\omega
_{2}^{3}-7y_{8}y_{7}^{2}-19y_{8}y_{7}z_{6}\omega _{2}-6y_{8}y_{7}y_{6}\omega
_{2}$

$\qquad +6y_{8}y_{7}y_{5}\omega _{2}^{2}-7y_{8}y_{7}y_{4}y_{3}\
-8y_{8}y_{7}y_{3}^{2}\omega _{2}+2y_{8}y_{7}y_{3}\omega _{2}^{4}+\
96y_{8}z_{6}^{2}\omega _{2}^{2}$

$\qquad +24y_{8}z_{6}y_{6}\omega _{2}^{2}+\
44y_{8}z_{6}y_{5}y_{3}-32y_{8}z_{6}y_{5}\omega
_{2}^{3}-59y_{8}z_{6}y_{4}^{2}+108y_{8}z_{6}y_{4}y_{3}\omega _{2}$

$\qquad -27y_{8}z_{6}y_{4}\omega _{2}^{4}-16y_{8}z_{6}y_{3}^{2}\omega
_{2}^{2}+2y_{8}y_{6}^{2}\omega
_{2}^{2}+6y_{8}y_{6}y_{5}y_{3}+y_{8}y_{6}y_{4}^{2}$

$\qquad -2y_{8}y_{6}y_{4}y_{3}\omega _{2}+y_{8}y_{6}y_{4}\omega _{2}^{4}+\
6y_{8}y_{5}y_{4}y_{3}\omega _{2}^{2}-3y_{8}y_{4}^{2}y_{3}^{2}+\
y_{8}y_{4}^{2}y_{3}\omega _{2}^{3}$

$\qquad -34y_{7}^{2}z_{6}\omega
_{2}^{2}+y_{7}z_{6}^{2}y_{3}-109y_{7}z_{6}^{2}\omega
_{2}^{3}-4y_{7}z_{6}y_{6}y_{3}+2y_{7}z_{6}y_{6}\omega _{2}^{3}$

$\qquad +8y_{7}z_{6}y_{5}y_{3}\omega _{2}-24y_{7}z_{6}y_{4}^{2}\omega
_{2}+4y_{7}z_{6}y_{4}y_{3}\omega _{2}^{2}+y_{7}y_{5}^{3}-51z_{6}^{3}y_{4}$

$\qquad -92z_{6}^{3}\omega
_{2}^{4}+102z_{6}^{2}y_{6}y_{4}-6z_{6}^{2}y_{6}y_{3}\omega
_{2}+8z_{6}^{2}y_{6}\omega _{2}^{4}+98z_{6}^{2}y_{5}y_{4}\omega _{2}$

$\qquad +96z_{6}^{2}y_{5}y_{3}\omega _{2}^{2}-153z_{6}^{2}y_{4}^{2}\omega
_{2}^{2}+55z_{6}^{2}y_{4}y_{3}^{2}-z_{6}^{2}y_{4}y_{3}\omega
_{2}^{3}-4z_{6}^{2}y_{3}^{3}\omega _{2}$

$\qquad +12z_{6}y_{6}^{2}y_{3}\omega _{2}-4z_{6}y_{6}^{2}\omega
_{2}^{4}-12z_{6}y_{6}y_{4}^{2}\omega
_{2}^{2}+8z_{6}y_{6}y_{4}y_{3}^{2}+2z_{6}y_{6}y_{4}y_{3}\omega _{2}^{3}$

$\qquad -2z_{6}y_{6}y_{3}^{3}\omega _{2}+y_{5}^{3}y_{4}\omega _{2}^{3})$,
\end{enumerate}

\noindent \textsl{and where }$z_{6}=2y_{6}+y_{3}^{2}+2\omega
_{2}y_{5}-3\omega _{2}^{2}y_{4}+2\omega _{2}^{3}y_{3}-\omega _{2}^{6}$.$%
\square $

\bigskip

We note that results in (3.5), (3.6) has been shown in \cite[Theorems 1,3]%
{DZ2}.

\section{Computing with Weyl invariants}

As mentioned earlier our approach to the ring $H^{\ast }(G/T)$ amounts to
apply Lemma 2.5 to the fibration (3.2). It requires in addition to Lemma 3.1
that

\begin{quote}
i) A presentation for the cohomology of the fiber space $P_{\{k\}}/T$;

ii) a set $\{\rho _{s}\}_{1\leq s\leq m_{1}}$ of relations on $H^{\ast
}(G/T) $ satisfying (2.4).
\end{quote}

\noindent These two tasks will be implemented in Lemmas 4.2 and 4.4,
respectively.

The Weyl group $W$ of a Lie group $G$ can be regarded as the subgroup of $%
Aut(H^{2}(G/T))$ generated by the elements $\sigma _{1},\cdots ,\sigma
_{n}\in Aut(H^{2}(G/T))$ whose action on set $\{\omega _{1},\cdots ,\omega
_{n}\}$ of weights is (see \cite[Section 2.1]{DZ2})

\begin{enumerate}
\item[(4.1)] $\sigma _{{i}}(\omega _{{k}})=\left\{ 
\begin{tabular}{l}
$\omega _{{i}}\text{ if }k\neq i\text{;}$ \\ 
$\omega _{{i}}-\sum\nolimits_{{1\leq j\leq n}}c_{{ij}}\omega _{{j}}\text{ if 
}k=i$,%
\end{tabular}%
\right. $ $1\leq i\leq n$,
\end{enumerate}

\noindent where $c_{{ij}}$ is the Cartan number relative to the pair $\beta
_{i},\beta _{j}$, $1\leq i,j\leq n$, of simple roots. Given a subgroup $%
W^{\prime }\subseteq W$ and a weight $\omega \in H^{2}(G/T)$ let $O(\omega
,W^{\prime })\subset H^{2}(G/T)$ be the $W^{\prime }$--orbit through $\omega 
$, and write $e_{r}(O(\omega ,W^{\prime }))\in $ $H^{\ast }(G/T)$ for the $%
r^{th}$ elementary symmetric functions on the set $O(\omega ,W^{\prime })$.
In the table below we define for each simple Lie group $G\neq Spin(n)$ a set 
$c_{r}(G)\in H^{\ast }(G/T)$ of polynomials in the weights $\omega
_{1},\cdots ,\omega _{n}$:

\begin{enumerate}
\item[(4.2)] {\footnotesize 
\begin{tabular}{l|l|l|l|l}
\hline\hline
$G$ & $SU(n),Sp(n)$ & $G_{2}$ & $F_{4}$ & $E_{n},n=6,7,8$ \\ \hline
$c_{r}(G)$ & $e_{r}(O(\omega _{1},W))$ & $e_{r}(O(\omega _{2},W_{\{1\}}))$ & 
$e_{r}(O(\omega _{4},W_{\{1\}}))$ & $e_{r}(O(\omega _{n},W_{\{2\}}))$ \\ 
\hline\hline
\end{tabular}%
},
\end{enumerate}

\noindent where $W_{\{i\}}$ is the Weyl group of the parabolic subgroup $%
P_{\{i\}}\subset G$ specified in (3.1).

\bigskip

\noindent \textbf{Example 4.1.} The expressions of $c_{r}(G)$ as polynomials
in the weights $\omega _{1},\cdots ,\omega _{n}$ can be concretely
presented. As examples, in the order of $G=SU(n),Sp(n),F_{4}$, $E_{6}$, we
get from the formula (4.1), together with the Cartan matrix of $G$ given in 
\cite[p.59]{H}, that

\begin{quote}
$O(\omega _{1},W){\footnotesize =}\left\{ {\footnotesize \omega }_{%
{\footnotesize 1}}{\footnotesize ,\omega }_{{\footnotesize k}}{\footnotesize %
-\omega }_{{\footnotesize k-1}}{\footnotesize ,-\omega }_{{\footnotesize n-1}%
}{\footnotesize \mid 2\leq k\leq n-1}\right\} $;

$O(\omega _{1},W){\footnotesize =}\left\{ {\footnotesize \pm \omega }_{%
{\footnotesize 1}}{\footnotesize ,\pm (\omega }_{{\footnotesize k}}%
{\footnotesize -\omega }_{{\footnotesize k-1}}{\footnotesize )\mid 2\leq
k\leq n}\right\} $;

$O(\omega _{4},W_{\{1\}}){\footnotesize =}\left\{ {\footnotesize \omega }%
_{_{4}}{\footnotesize ,\omega }_{{\footnotesize 3}}{\footnotesize -\omega }_{%
{\footnotesize 4}}{\footnotesize ,\omega }_{{\footnotesize 2}}{\footnotesize %
-\omega }_{{\footnotesize 3}}{\footnotesize ,\omega }_{{\footnotesize 1}}%
{\footnotesize -\omega }_{{\footnotesize 2}}{\footnotesize +\omega }_{%
{\footnotesize 3}}{\footnotesize ,\omega }_{{\footnotesize 1}}{\footnotesize %
-\omega }_{{\footnotesize 3}}{\footnotesize +\omega }_{{\footnotesize 4}}%
{\footnotesize ,\omega }_{{\footnotesize 1}}{\footnotesize -\omega }_{%
{\footnotesize 4}}\right\} $;

$O(\omega _{6},W_{\{2\}}){\footnotesize =}\left\{ {\footnotesize \omega }_{%
{\footnotesize 6}}{\footnotesize ,\omega }_{{\footnotesize 5}}{\footnotesize %
-\omega }_{{\footnotesize 6}}{\footnotesize ,\omega }_{{\footnotesize 4}}%
{\footnotesize -\omega }_{{\footnotesize 5}}{\footnotesize ,\omega }_{%
{\footnotesize 2}}{\footnotesize +\omega }_{{\footnotesize 3}}{\footnotesize %
-\omega }_{{\footnotesize 4}}{\footnotesize ,\omega }_{{\footnotesize 1}}%
{\footnotesize +\omega }_{{\footnotesize 2}}{\footnotesize -\omega }_{%
{\footnotesize 3}}{\footnotesize ,\omega }_{{\footnotesize 2}}{\footnotesize %
-\omega }_{{\footnotesize 1}}\right\} $.$\square $
\end{quote}

In the following results we clarify the roles of the polynomials $c_{r}(G)$.

\bigskip

\noindent \textbf{Lemma 4.2. }\textsl{If }$G=SU(n)$\textsl{\ or }$Sp(n)$%
\textsl{\ the inclusion }$\omega _{i}\subset H^{2}(G/T)$ \textsl{induces
ring isomorphisms}

\begin{enumerate}
\item[(4.2)] $H^{\ast }({\small SU(n)/T})=\mathbb{Z}\left[ \omega
_{1},\cdots ,\omega _{n-1}\right] /\left\langle c_{2},\cdots
,c_{n}\right\rangle $\textsl{,} $c_{r}=c_{r}({\small SU(n)})$,

\item[(4.3)] $H^{\ast }(Sp(n)/T)=\mathbb{Z}\left[ \omega _{1},\cdots ,\omega
_{n}\right] /\left\langle c_{2},\cdots ,c_{2n}\right\rangle $\textsl{,} $%
c_{2r}=c_{2r}(Sp(n))$.
\end{enumerate}

\noindent \textbf{Proof.} For $G=SU(n)$ or $Sp(n)$ we have by Borel \cite{B}
that

\begin{quote}
$H^{\ast }(G/T)=\mathbb{Z}\left[ \omega _{1},\cdots ,\omega _{n}\right]
/\left\langle \mathbb{Z}\left[ \omega _{1},\cdots ,\omega _{n}\right]
^{+,W}\right\rangle $,
\end{quote}

\noindent where $\mathbb{Z}\left[ \omega _{1},\cdots ,\omega _{n}\right]
^{+,W}$\ denotes the set of $W$--invariants in positive degrees. The lemma
is verified by the classical results that the sets of polynomials $c_{r}(G)$
in (4.2) and (4.3) generates the subrings $\mathbb{Z}\left[ \omega
_{1},\cdots ,\omega _{n}\right] ^{+,W}$.$\square $

\bigskip

For an exceptional Lie group $G$ let $P_{\{k\}}\subset G$ be the parabolic
subgroup by (3.1), and consider the corresponding fibration

\begin{quote}
$P_{\{k\}}^{s}/T^{\prime }\overset{i}{\hookrightarrow }G/T\overset{\pi }{%
\rightarrow }G/P_{\{k\}}$
\end{quote}

\noindent in flag manifolds, where $P_{\{k\}}^{s}$\ is the simple part of
the group $P_{\{k\}}$, and where $T^{\prime }$ is the maximal torus on $%
P_{\{k\}}^{s}$ corresponding to $T$.

\bigskip

\noindent \textbf{Lemma 4.3. }\textsl{For each exceptional Lie group }$G$%
\textbf{\ }\textsl{the polynomials }$c_{r}(G)\in H^{\ast }(G/T)$\textsl{\
defined in table (4.2) satisfy the following relations}

\textsl{i)} $c_{r}(G)\in \func{Im}[\pi ^{\ast }:H^{\ast
}(G/P_{\{k\}})\rightarrow H^{\ast }(G/T)]$\textsl{;}

\textsl{ii) }$i^{\ast }c_{r}(G)=c_{r}(P_{\{k\}}^{s}).$

\bigskip

\noindent \textbf{Proof.} For any parabolic subgroup $P\subset G$ with Weyl
group $W(P)$ the induced map $\pi ^{\ast }$ is injective by Lemma 2.3, and
satisfies that

\begin{enumerate}
\item[(4.4)] $\func{Im}\pi ^{\ast }=H^{\ast }(G/T)^{W(P)}$ (see \cite[%
Proposition 5.1]{BGG}),
\end{enumerate}

\noindent where $H^{\ast }(G/T)^{W(P)}\subset H^{\ast }(G/T)$ is the subring
of $W(P)$--invariants. Property i) follows from $c_{r}(G)\in H^{\ast
}(G/T)^{W_{\{k\}}}$ by Definition (4.2), where $k=1$ for $G=G_{2}$ or $F_{4}$%
, and $k=2$ for $G=E_{n}$ with $n=6,7,8$.

By the $W_{\{k\}}$--equivariance of the induced map $i^{\ast }$ relation ii)
is verified by $i^{\ast }(O(\omega _{t},W_{\{k\}}))=O(i^{\ast }\omega
_{t},W_{\{k\}})$, where $t=2,4,6,7,8$ in accordance to $%
G=G_{2},F_{4},E_{6}.E_{7}$ and $E_{8}$.$\square $

\bigskip

\noindent \textbf{Lemma 4.4.} \textsl{For each exceptional Lie group }$G$%
\textsl{\ a set }$\{\rho _{s}\}$ \textsl{of relations on the ring }$H^{\ast
}(G/T)$\textsl{\ satisfying the property (2.4) are given in the table below:}

\begin{center}
\begin{tabular}{l|l}
\hline\hline
$G$ & $\{\rho _{s}\}$ \\ \hline
$G_{2}$ & $3\omega _{1}-c_{1}$; $3\omega _{1}^{2}-c_{1}$ \\ \hline
$F_{4}$ & $3\omega _{1}-c_{1}$; $4\omega _{1}^{2}-c_{2}$; $6y_{3}-c_{3}$; $%
3y_{4}+2\omega _{1}y_{3}-c_{4}$; $\omega _{1}y_{4}-c_{5}$; $y_{6}-c_{6}$ \\ 
\hline
$E_{6}$ & 
\begin{tabular}{l}
$3\omega _{2}-c_{1}$; $4\omega _{2}^{2}-c_{2}$; $2y_{3}+2\omega
_{2}^{3}-c_{3}$; $3y_{4}+\omega _{2}^{4}-c_{4}$; \\ 
$3\omega _{2}y_{4}-2\omega _{2}^{2}y_{3}+\omega _{2}^{5}-c_{5}$; $%
y_{6}-c_{6} $%
\end{tabular}
\\ \hline
$E_{7}$ & 
\begin{tabular}{l}
$3\omega _{2}-c_{1}$; $4\omega _{2}^{2}-c_{2}$; $2y_{3}+2\omega
_{2}^{3}-c_{3}$; $3y_{4}+\omega _{2}^{4}-c_{4}$; \\ 
$2y_{5}+3\omega _{2}y_{4}-2\omega _{2}^{2}y_{3}+\omega _{2}^{5}-c_{5}$; $%
y_{6}+2\omega _{2}y_{5}-c_{6}$; $y_{7}-c_{7}$%
\end{tabular}
\\ \hline
$E_{8}$ & 
\begin{tabular}{l}
$3\omega _{2}-c_{1}$; $4\omega _{2}^{2}-c_{2}$; $2y_{3}+2\omega
_{2}^{3}-c_{3}$; $3y_{4}+\omega _{2}^{4}-c_{4}$; \\ 
$2y_{5}+3\omega _{2}y_{4}\ -2\omega _{2}^{2}y_{3}+\omega _{2}^{5}-c_{5}$; \\ 
$5{y_{6}+2{y_{3}^{2}+6{\omega }_{2}{y_{5}-6\omega }_{2}^{2}{y_{4}+4\omega }%
_{2}^{3}{y_{3}}}-2{\omega }_{2}^{6}-}c_{6}$; \\ 
$y_{7}+4\omega _{2}y_{6}+2\omega _{2}y_{3}^{2}\ +4\omega
_{2}^{2}y_{5}-6\omega _{2}^{3}y_{4}+4\omega _{2}^{4}y_{3}-2\omega
_{2}^{7}-c_{7}$; \\ 
$y_{8}-c_{8}$%
\end{tabular}
\\ \hline\hline
\end{tabular}
\end{center}

\noindent \textsl{where the }$y_{i}$\textsl{'s are the }$\pi ^{\ast }$%
\textsl{--images of the} \textsl{Schubert classes on the }$G/P_{\{k\}}$%
\textsl{\ specified by (3.3)}, $c_{r}=c_{r}(G),$ \textsl{and where the set }$%
\{\rho _{s}\}$\textsl{\ are presented by the order of the degrees of the
enclosed polynomials }$\rho _{s}$\textsl{.}

\bigskip

\noindent \textbf{Proof.} By i) of Lemma 4.3 and by the injectivity of the
map $\pi ^{\ast }$ we can regard $c_{r}(G)\in H^{\ast }(G/P_{\{k\}})$, see
Convention 2.4. Moreover since $c_{r}(G)$ is a polynomial in the Schubert
classes $\omega _{1},\cdots ,\omega _{n}$ the\textsl{\ }package of\textsl{\
"Giambelli polynomials"} \cite[Section 2.6]{DZ2}\textsl{\ }is functional to
expand it as a polynomial in the special Schubert classes on $H^{\ast
}(G/P_{\{k\}})$ given in (3.3). This yields the relations $\rho _{r}$'s on
the ring $H^{\ast }(G/T)$ presented in table.

Finally, by Lemma 4.2 and ii) of Lemma 4.3 property (2.4) is satisfied by
the set $\left\{ \rho _{r}\right\} $ of relations on $H^{\ast }(G/T)$.$%
\square $

\bigskip

\noindent \textbf{Remark 4.5.} Results in Lemma 4.4 has geometric
interpretations. Taking $G=E_{n}$ with $n=6,7,8$ as examples the subgroup $%
P_{(2)}=SU(n)\cdot S^{1}$ has a canonical $n$--dimensional complex
representation that gives rise to a complex $n$--bundle $\xi _{n}$ on the
Grassmannian $E_{n}/P_{\{2\}}$ \cite{AH}. It can be shown that if we let $%
c_{r}(\xi _{n})\in H^{r}(E_{n}/P_{\{2\}})$ be the $r^{th}$ Chern class of $%
\xi _{n}$, $1\leq r\leq n$, then

\begin{quote}
$c_{r}(\xi _{n})=c_{r}(G)$, $1\leq r\leq n$.
\end{quote}

\noindent In this regard the relations $\rho _{s}=0$ indicate formulae that
express the Chern classes $c_{r}(\xi _{n})$ by the special Schubert classes
on $E_{n}/P_{\{2\}}$.

Note that if we let $p:\mathbb{C}P(\xi _{n})\rightarrow E_{n}/P_{\{2\}}$ be
the complex projective bundle associated to $\xi _{n}$, then $\mathbb{C}%
P(\xi _{n})=E_{n}/P_{\{2,n\}}$, and the projection $p$ agrees with the
bundle map induced by the inclusion $P_{\{2,n\}}\subset P_{\{2\}}\subset
E_{n}$ of parabolic subgroups.$\square $

\section{The ring $H^{\ast }(G/T)$ for exceptional Lie groups}

Summarizing the computation of Sections 3 and 4 we have associated each
exceptional Lie group $G$ with a fibration $P_{\{k\}}^{s}/T^{\prime
}\hookrightarrow G/T\rightarrow G/P_{\{k\}}$ in which presentations of the
cohomologies of the base and fiber spaces by Schubert classes have been
obtained in Lemma 3.1 and Lemma 4.2, respectively. In addition, a set of
relations on $H^{\ast }(G/T)$ satisfying the condition (2.4) has been
determined in Lemma 4.4. Therefore, Lemma 2.5 is directly applicable to
yield the following result, where the $y_{i}$'s are the Schubert classes on $%
G/P_{\{k\}}$ given by (3.3), and where $c_{r}=c_{r}(G)$ as in Lemma 4.4.

\bigskip

\noindent \textbf{Theorem 5.1. }\textsl{For each exceptional Lie group }$G$ 
\textsl{the cohomology ring }$H^{\ast }(G/T)$\textsl{\ has the following
presentation}

\begin{enumerate}
\item[(5.1)] $H^{\ast }(G_{2}/T)=\mathbb{Z}[\omega _{1},\omega
_{2},y_{3}]/\left\langle \rho _{2},r_{3},r_{6}\right\rangle $\textsl{, where}

$\rho _{2}=3\omega _{1}^{2}-3\omega _{1}\omega _{2}+\omega _{2}^{2}$;

$r_{3}=2y_{3}-\omega _{1}^{3}$;$\quad $

$r_{6}=y_{3}^{2}$.

\item[(5.2)] $H^{\ast }(F_{4}/T)=\mathbb{Z}[\omega _{1},\omega _{2},\omega
_{3},\omega _{4},y_{3},y_{4}]/\left\langle \rho _{2},\rho
_{4},r_{3},r_{6},r_{8},r_{12}\right\rangle $\textsl{, where}

$\rho _{2}=c_{2}-4\omega _{1}^{2}$;

$\rho _{4}=3y_{4}+2\omega _{1}y_{3}-c_{4}$;

$r_{3}=2y_{3}-\omega _{1}^{3}$;

$r_{6}=y_{3}^{2}+2c_{6}-3\omega _{1}^{2}y_{4}$;

$r_{8}=3y_{4}^{2}-\omega _{1}^{2}c_{6}$;

$r_{12}=y_{4}^{3}-c_{6}^{2}$.

\item[(5.3)] $H^{\ast }(E_{6}/T)=\mathbb{Z}[\omega _{1},\cdots ,\omega
_{6},y_{3},y_{4}]/\left\langle \rho _{2},\rho _{3},\rho _{4},\rho
_{5},r_{6},r_{8},r_{9},r_{12}\right\rangle $\textsl{, where}

$\rho _{2}=4\omega _{2}^{2}-c_{2}$;

$\rho _{3}=2y_{3}+2\omega _{2}^{3}-c_{3}$;

$\rho _{4}=3y_{4}+\omega _{2}^{4}-c_{4}$;

$\rho _{5}=2\omega _{2}^{2}y_{3}-\omega _{2}c_{4}+c_{5}$;

$r_{6}=y_{3}^{2}-\omega _{2}c_{5}+2c_{6}$;

$r_{8}=3y_{4}^{2}-2c_{5}y_{3}-\omega _{2}^{2}c_{6}+\omega _{2}^{3}c_{5};$

$r_{9}=2y_{3}c_{6}-\omega _{2}^{3}c_{6}$;

$r_{12}=y_{4}^{3}-c_{6}^{2}$.

\item[(5.4)] $H^{\ast }(E_{7}/T)=\mathbb{Z}[\omega _{1},\cdots ,\omega
_{7},y_{3},y_{4},y_{5},y_{9}]/\left\langle \rho _{2},\rho _{3},\rho
_{4},\rho _{5},r_{i}\right\rangle $\textsl{, where }$i\in \{6,$ $%
8,9,10,12,14,18\}$\textsl{, and where}

$\rho _{2}=4\omega _{2}^{2}-c_{2}$;

$\rho _{3}=2y_{3}+2\omega _{2}^{3}-c_{3}$;

$\rho _{4}=3y_{4}+\omega _{2}^{4}-c_{4}$;

$\rho _{5}=2y_{5}-2\omega _{2}^{2}y_{3}+\omega _{2}c_{4}-c_{5}$;

$r_{6}=y_{3}^{2}-\omega _{2}c_{5}+2c_{6}$;

$r_{8}=3y_{4}^{2}+2y_{3}y_{5}-2y_{3}c_{5}+2\omega _{2}c_{7}-\omega
_{2}^{2}c_{6}+\omega _{2}^{3}c_{5}$;

$r_{9}=2{y_{9}}+2{y_{4}y_{5}}-2{y_{3}c_{6}}-{\omega _{2}^{2}c_{7}}+{\omega
_{2}^{3}c_{6}}$;

$r_{10}=y_{5}^{2}-2y_{3}c_{7}+\omega _{2}^{3}c_{7}$;

$r_{12}=y_{4}^{3}-4y_{5}c_{7}-c_{6}^{2}-2y_{3}y_{9}-2y_{3}y_{4}y_{5}+2\omega
_{2}y_{5}c_{6}+3\omega _{2}y_{4}c_{7}+c_{5}c_{7}$;

$r_{14}=c_{7}^{2}-2y_{5}y_{9}+2y_{3}y_{4}c_{7}-\omega _{2}^{3}y_{4}c_{7}$;

$%
r_{18}=y_{9}^{2}+2y_{5}c_{6}c_{7}-y_{4}c_{7}^{2}-2y_{4}y_{5}y_{9}+2y_{3}y_{5}^{3}-5\omega _{2}y_{5}^{2}c_{7} 
$.

\item[(5.5)] $H^{\ast }(E_{8}/T)=\mathbb{Z}[\omega _{1},\cdots ,\omega
_{8},y_{3},y_{4},y_{5},y_{6},y_{9},y_{10},y_{15}]/\left\langle \rho
_{i},r_{j}\right\rangle $\textsl{, where }$i\in \{2,3,4,5,6\}$\textsl{, }$%
j\in \{8,9,10,12,14,15,18,20,24,30\}$\textsl{\ and where}

$\rho _{2}=4\omega _{2}^{2}-c_{2}$;

$\rho _{3}=2y_{3}+2\omega _{2}^{3}-c_{3}$;

$\rho _{4}=3y_{4}+\omega _{2}^{4}-c_{4}$;

$\rho _{5}=2y_{5}-2\omega _{2}^{2}y_{3}+\omega _{2}c_{4}-c_{5}$;

$\rho _{6}=5y_{6}+2y_{3}^{2}+10\omega _{2}y_{5}-2\omega _{2}c_{5}-c_{6}$;

$r_{8}=3{c_{8}}-3{y_{4}^{2}}-2{y_{3}y_{5}}+2{y_{3}}c_{5}-{2\omega _{2}c_{7}}+%
{\omega _{2}^{2}c_{6}-\omega _{2}^{3}c_{5}}$;

$r_{9}=2{y_{9}}+2{y_{4}y_{5}}-2{y_{3}y_{6}}-4{{\omega _{2}}y_{3}y_{5}+{%
\omega _{2}c_{8}}-\omega _{2}^{2}c_{7}}+{\omega _{2}^{3}c_{6}}$;

$r_{10}=3{{y_{10}}+6{y_{4}z_{6}}}-2{y_{5}^{2}}-2{y_{3}c_{7}}-{\omega
_{2}^{2}c_{8}}+{\omega _{2}^{3}c_{7}}$;

$r_{12}={y_{4}^{3}}-2{y_{3}y_{9}}-2{y}_{{3}}^{2}{z_{6}}-2{y_{3}y_{4}y_{5}}%
-c_{6}^{2}+4{c_{6}z_{6}}\ {+2\omega _{2}y_{5}}c_{6}-2{\omega _{2}c_{5}z_{6}}$

$\qquad -4{y_{5}c_{7}}+3{\omega _{2}y_{4}c_{7}}+c_{5}{c_{7}}-6{y_{4}c_{8}}$;

$r_{14}\equiv c_{{7}}^{2}+3y_{{4}}{y_{10}}+6{y_{4}^{2}z_{6}}-2y_{{3}}y_{{5}%
}c_{{6}}-{\omega _{2}^{2}}y_{{5}}c_{{7}}+{\omega _{2}^{3}}y_{{5}}c_{{6}}%
\func{mod}c_{8}$;

$r_{15}\equiv 2y_{15}+2y_{5}({y_{10}-y_{5}^{2}}+2{y_{4}z_{6}}%
)+2y_{3}z_{6}^{2}-2\omega _{2}^{3}z_{6}^{2}\func{mod}c_{8}$;

$r_{18}\equiv
y_{9}^{2}-10z_{6}^{3}+5y_{5}z_{6}y_{7}-y_{4}y_{7}^{2}-2y_{4}y_{5}y_{9}+4y_{3}z_{6}y_{9}+6y_{3}y_{5}y_{10} 
$

$\qquad -4y_{3}y_{5}^{3}\
+8y_{3}y_{4}y_{5}z_{6}+4y_{3}^{2}z_{6}^{2}-y_{3}^{2}y_{5}y_{7}-3\omega
_{2}y_{7}y_{10}+10\omega _{2}y_{5}z_{6}^{2}$

$\qquad -10\omega _{2}y_{4}z_{6}y_{7}-2\omega
_{2}y_{3}^{2}y_{5}z_{6}-6\omega _{2}^{2}z_{6}y_{10}-16\omega
_{2}^{2}y_{4}z_{6}^{2}+3\omega _{2}^{2}y_{4}y_{5}y_{7}$

$\qquad -3\omega _{2}^{3}y_{5}y_{10}+2\omega _{2}^{3}y_{5}^{3}\func{mod}%
c_{8} $;

$r_{20}\equiv 3({y_{5}^{2}-y_{10}}-2{y_{4}z_{6}})^{2}-{y_{5}^{3}\rho }_{5}%
\func{mod}c_{8};$

$r_{24}\equiv 5(2{y_{6}}+{y_{3}^{2}}+4\omega _{2}{y_{5}}-\omega
_{2}c_{5})^{4}-2{y_{3}^{7}\rho }_{3}\func{mod}c_{8};$

$r_{30}\equiv
(-y_{15}-y_{5}y_{10}+y_{5}^{3}-2y_{4}y_{5}z_{6}-y_{3}z_{6}^{2}+\omega
_{2}^{3}z_{6}^{2})^{2}+({y_{5}^{2}-y_{10}}-2{y_{4}z_{6}})^{3}$

$\qquad +{(2{y_{6}}+{y_{3}^{2}}+4\omega _{2}{y_{5}}-\omega _{2}c_{5})}^{5}-6{%
{y_{6}^{4}}\rho }_{6}\func{mod}c_{8},$
\end{enumerate}

\noindent \textsl{in which} $z_{6}=2{y_{6}}+{y_{3}^{2}}+4\omega _{2}{y_{5}}%
-\omega _{2}c_{5}$.$\square $

\bigskip

Concerning the formulation of the presentations (5.1)--(5.5) in Theorem 5.1
we remark that

a) certain Schubert classes $y_{k}$ on the base space $G/P_{\{k\}}$ can be
eliminated against appropriate relations of the type $\rho _{k}$. As example
when $G=E_{7}$ the generators $y_{6},y_{7}$ and the relations $\rho
_{6},\rho _{7}$ can be excluded by the formulae of $\rho _{6}$ ($%
y_{6}=c_{6}-2\omega _{2}y_{5}$) and $\rho _{7}$ ($y_{7}=c_{7}$) in Lemma 4.4.

b) for simplicity the relations $r_{k}$ on the ring $H^{\ast }(E_{8}/T)$
with $k\geq 14$ are presented after module $c_{8}$, \ while their full
expressions have been recorded in (3.7);

c) without altering the ideal, higher degree relations of the type $r_{i}$
may be simplified using the lower degree relations. The main idea of
performing such simplifications is the following one: for two ordered
subsets $\{f_{i}\}_{1\leq i\leq n}$ and $\{h_{i}\}_{1\leq i\leq n}$ of a
graded polynomial ring with

\begin{quote}
$\deg f_{1}<\cdots <\deg f_{n}$ and $\deg h_{1}<\cdots <\deg h_{n}$
\end{quote}

\noindent write $\{h_{i}\}_{1\leq i\leq n}\thicksim \{f_{i}\}_{1\leq i\leq
n} $ to denote the statements that $\deg h_{i}=\deg f_{i}$ and that $%
(f_{i}-h_{i})\in \left\langle f_{j}\right\rangle _{1\leq j<i}$. Then

\begin{enumerate}
\item[(5.6)] \noindent $\{f_{i}\}_{1\leq i\leq n}$ $\sim $ $\{h_{i}\}_{1\leq
i\leq n}$ implies that $\left\langle h_{1},\cdots ,h_{n}\right\rangle
=\left\langle f_{1},\cdots ,f_{n}\right\rangle $.$\square $
\end{enumerate}

\section{Proofs of Theorems 1.2 and 1.3}

\noindent \textbf{Proof of Theorem 1.2.} If $G=SU(n)$ or $Sp(n)$ we have $%
m=0 $, and the presentation (1.2) is shown by Lemma 4.2. If $%
G=G_{2},F_{4},E_{6}$ and $E_{7}$ the formula (1.2) is verified by
(5.1)--(5.4), in the light of the relation (5.6).

For $G=E_{8}$ the presentation (5.5) can be summarized as

\begin{enumerate}
\item[(6.1)] $H^{\ast }(E_{8}/T)=\mathbb{Z}[\omega _{1},\cdots ,\omega
_{8},y_{r}]/\left\langle e_{i},f_{j},g_{t},\phi \right\rangle _{1\leq i\leq
3;1\leq j\leq 7,t=1,2,3,5}$
\end{enumerate}

\noindent where\textsl{\ }$r\in \{3,4,5,6,9,10,15\}$ and where

\begin{quote}
i) $e_{i}\in \left\langle \omega _{1},\cdots ,\omega _{8}\right\rangle $, $%
1\leq i\leq 3$;

ii) $f_{j}$\ $=$\ $p_{j}y_{d_{j}}+\alpha _{j}$, $p_{j}\in \{2,3,5\}$, $%
\alpha _{j}\in \left\langle \omega _{1},\cdots ,\omega _{8}\right\rangle $, $%
1\leq j\leq 7$;

iii) $g_{t}=y_{d_{t}}^{k_{t}}+\beta _{t}$ with $\beta _{t}\in \left\langle
\omega _{1},\cdots ,\omega _{8}\right\rangle $, $t=1,2,3,5$;

iv) $\phi =2y_{6}^{5}-y_{10}^{3}+y_{15}^{2}+\beta $\ with $\beta \in
\left\langle \omega _{1},\cdots ,\omega _{8}\right\rangle $.
\end{quote}

\noindent Comparing (1.2) with (6.1) we find that

\begin{quote}
a) the polynomial $\phi $ in iv) does not belong to any of the three types $%
e_{i}$, $f_{j}$, $g_{j}$ of relations in Theorem 1.2;

b) the polynomials $g_{4},g_{6},g_{7}$ in (1.2) required to couple $%
f_{4},f_{6},f_{7}$ (see Theorem 1.2) are absent in iii).
\end{quote}

\noindent However, if we set

\begin{enumerate}
\item[(6.2)] $\left\{ 
\begin{tabular}{l}
$g_{4}=-12\phi +5y_{6}^{4}f_{4}-4y_{10}^{2}f_{6}+6y_{15}f_{7}$ \\ 
$g_{6}=-10\phi +4y_{6}^{4}f_{4}-3y_{10}^{2}f_{6}+5y_{15}f_{7}$ \\ 
$g_{7}=15\phi -6y_{6}^{4}f_{4}+5y_{10}^{2}f_{6}-7y_{15}f_{7}$%
\end{tabular}%
\ \ \right. $.
\end{enumerate}

\noindent then the obvious properties

\begin{quote}
$g_{4},g_{6},g_{7}\in \left\langle e_{i},f_{k},g_{s},\phi \right\rangle $; $%
\phi =2g_{4}-g_{6}+g_{7}\in \left\langle e_{i};f_{j},g_{j}\right\rangle
_{1\leq i\leq 3,1\leq j\leq 7}$,
\end{quote}

\noindent with $1\leq i\leq 3$; $1\leq k\leq 7$, $s=1,2,3,5$, imply the
relation

\begin{quote}
$\left\langle e_{i};f_{j},g_{s},\phi \right\rangle _{1\leq i\leq 3;1\leq
j\leq 7,s=1,2,3,5}=\left\langle e_{i};f_{j},g_{j}\right\rangle _{1\leq i\leq
3,1\leq j\leq 7}$.
\end{quote}

\noindent That is, the formula (1.2) for $G=E_{8}$ is indeed identical to
(6.1).

For the remaining case $G=Spin(m)$ let $y_{k}$ be the Schubert class on $%
Spin(2n)/T$ associated to the element $w_{k}=\sigma \lbrack n-k,\cdots
,n-2,n-1]$ in the Weyl group of $Spin(2n)$, $2\leq k\leq n-1$. According to
Marlin \cite[Proposition 3]{M} one has the presentation

\begin{quote}
$H^{\ast }(Spin(2n)/T)=\mathbb{Z}[\omega _{1},\cdots ,\omega
_{n},y_{2},\cdots ,y_{n-1}]/\left\langle \delta _{i},\xi _{j},\mu
_{k}\right\rangle $
\end{quote}

\noindent with

\begin{quote}
$\delta _{i}:=2y_{i}-c_{i}(\omega _{1},\cdots ,\omega _{n})$, $1\leq i\leq
n-1$,

$\xi _{j}:=y_{2j}+(-1)^{j}y_{j}^{2}+2\sum\limits_{1\leq r\leq
j-1}(-1)^{r}y_{r}y_{2j-r}$, $1\leq j\leq \left[ \frac{n-1}{2}\right] $,

$\mu _{k}:=(-1)^{k}y_{k}^{2}+2\sum\limits_{2k-n+1\leq r\leq
k-1}(-1)^{r}y_{r}y_{2k-r}$, $\left[ \frac{n}{2}\right] \leq k\leq n-1$,
\end{quote}

\noindent where $c_{i}(\omega _{1},\cdots ,\omega _{n})$ is the $i^{th}$
elementary symmetric function on the orbit set

\begin{quote}
$O{\footnotesize (\omega }_{n},{\small W}{\footnotesize )=}\left\{ 
{\footnotesize \omega }_{n}{\footnotesize ,\omega }_{i}{\footnotesize %
-\omega }_{i-1}{\footnotesize ,\omega }_{n-1}{\footnotesize +\omega }_{n}%
{\footnotesize -\omega }_{n-2}{\footnotesize ,\omega }_{n-1}{\footnotesize %
-\omega }_{n}{\footnotesize ,2\leq i\leq n-2}\right\} $,
\end{quote}

\noindent In view of the relations of the type $\xi _{j}$ we note that the
generators $y_{2j}$ with $1\leq j\leq \left[ \frac{n-1}{2}\right] $ can be
eliminated to yield the compact presentation

\begin{enumerate}
\item[(6.3)] $H^{\ast }(Spin(2n)/T)=\mathbb{Z}[\omega _{1},\cdots ,\omega
_{n},y_{3},y_{5},\cdots ,y_{2\left[ \frac{n-1}{2}\right] -1}]/\left\langle
\delta _{i}^{\prime },\mu _{k}^{\prime }\right\rangle $,
\end{enumerate}

\noindent where $\delta _{i}^{\prime }$ and $\mu _{k}^{\prime }$ are the
polynomials obtained from $\delta _{i}$ and $\mu _{k}$ by replacing all the
classes $y_{2r}$ by the polynomials

\begin{quote}
$(-1)^{r-1}y_{r}^{2}+2\sum\limits_{1\leq k\leq r-1}(-1)^{k-1}y_{k}y_{2r-k}$
(by the relation $\xi _{r}$).
\end{quote}

\noindent For $G=Spin(2n)$ formula (1.2) and is verified by (6.3).
Similarly, one obtains formula (1.2) for the group $G=Spin(2n+1)$ from \cite[%
Proposition 2]{M}.

It remains to show that the numbers $n$ and $m$ in (1.2) satisfy the
constraint $h(G,T)=n+m+1$. This will be done in the proof of Theorem 1.3.$%
\square $

\bigskip

The sets of integers appearing in the formula (1.2)

\begin{enumerate}
\item[(6.4)] $\{k,m\}$, $\{\deg e_{i}\}_{1\leq i\leq k}$, $%
\{d_{j};p_{j};k_{j}\}_{1\leq j\leq m}$,
\end{enumerate}

\noindent can be shown to be invariants of the corresponding Lie group $G$,
and will be called \textsl{the basic data} of $G$. With the formula (1.2)
being made explicit all the simple Lie groups in Lemma 4.2, formulae
(5.1)--(5.5), and in (6.3) one gets that

\bigskip

\noindent \textbf{\noindent Corollary 6.1. }\textsl{The basic data of the }$%
1 $\textsl{--connected simple Lie groups are}

\begin{center}
{\footnotesize 
\begin{tabular}{l|llll}
\hline\hline
$G$ & $SU(n+1)$ & $Sp(n)$ & $Spin(2n)$ & $Spin(2n+1)$ \\ \hline
$\{k,m\}$ & $\{n,0\}$ & $\{n,0\}$ & $\{[\frac{n+3}{2}],[\frac{n-2}{2}]\}$ & $%
\{[\frac{n+2}{2}],[\frac{n-1}{2}]\}$ \\ 
$\{\deg e_{i}\}$ & $\{2i+2\}$ & $\{4i\}$ & $\{4t,2n,2^{[\log
_{2}(n-1)]+2}\}_{1\leq t\leq \lbrack \frac{n-1}{2}]}$ & $\{4t,2^{[\log
_{2}n]+2}\}_{1\leq t\leq \lbrack \frac{n}{2}]}$ \\ 
$\{d_{j}\}$ &  &  & $\{4j+2\}$ & $\{4j+2\}$ \\ 
$\{p_{j}\}$ &  &  & $\{2,\cdots ,2\}$ & $\{2,\cdots ,2\}$ \\ 
$\{k_{j}\}$ &  &  & $\{2^{[\log _{2}\frac{n-1}{2j+1}]+1}\}$ & $\{2^{[\log
_{2}\frac{n}{2j+1}]+1}\}$ \\ \hline\hline
\end{tabular}
}

{\small Table 1. Basic data for the classical groups.}

{\footnotesize 
\begin{tabular}{l}
\hline\hline
\begin{tabular}{l|lllll}
$G$ & $G_{2}$ & $F_{4}$ & $E_{6}$ & $E_{7}$ & $E_{8}$ \\ \hline
$\{k,m\}$ & $\{1,1\}$ & $\{2,2\}$ & $\{4,2\}$ & $\{3,4\}$ & $\{3,7\}$ \\ 
$\{\deg e_{i}\}$ & $\{4\}$ & $\{4,16\}$ & $\{4,10,16,18\}$ & $\{4,16,28\}$ & 
$\{4,16,28\}$ \\ 
$\{d_{j}\}$ & $\{6\}$ & $\{6,8\}$ & $\{6,8\}$ & $\{6,8,10,18\}$ & $%
\{6,8,10,12,18,20,30\}$ \\ 
$\{p_{j}\}$ & $\{2\}$ & $\{2,3\}$ & $\{2,3\}$ & $\{2,3,2,2\}$ & $%
\{2,3,2,5,2,3,2\}$ \\ 
$\{k_{j}\}$ & $\{2\}$ & $\{2,3\}$ & $\{2,3\}$ & $\{2,3,2,2\}$ & $%
\{8,3,4,5,2,3,2\}$%
\end{tabular}
\\ \hline\hline
\end{tabular}%
}

{\small Table 2. Basic data for exceptional Lie groups}{\footnotesize .}$%
\square $
\end{center}

For a Lie group $G$ we set $\mathcal{A}(G):=H^{\ast }(G/T)/\left\langle
\omega _{1},\cdots ,\omega _{n}\right\rangle $. In the associated short
exact sequence of graded rings

\begin{enumerate}
\item[(6.4)] $0\rightarrow \left\langle \omega _{1},\cdots ,\omega
_{n}\right\rangle \rightarrow H^{\ast }(G/T)\overset{p}{\rightarrow }%
\mathcal{A}(G)\rightarrow 0$.
\end{enumerate}

\noindent the quotient map $p$ is clearly given by $p(\alpha )=\alpha \mid
_{\omega _{1}=\cdots =\omega _{n}=0}$, $\alpha \in H^{\ast }(G/T)$. It
follows from the formulae (1.2) and (1.3) that

\begin{enumerate}
\item[(6.5)] $\mathcal{A}(G)=\left\{ 
\begin{tabular}{l}
$\frac{\mathbb{Z}[y_{d_{1}},\cdots ,y_{d_{m}}]}{\left\langle p_{i}\cdot
y_{d_{i}},y_{d_{i}}^{k_{i}}\right\rangle _{1\leq i\leq m}}$ \quad if $G\neq
E_{8}$ \\ 
$\frac{\mathbb{Z}[y_{d_{1}},\cdots ,y_{d_{7}}]}{\left\langle
p_{i}y_{d_{i}},y_{d_{t}}^{k_{t}},2y_{d_{4}}^{5}-y_{d_{6}}^{3}+y_{d_{7}}^{2}%
\right\rangle _{1\leq i\leq 7,t=1,2,3,5}}\quad $if $G=E_{8}$.%
\end{tabular}%
\right. $
\end{enumerate}

\noindent Inputting the values of the data $\{d_{j};p_{j};k_{j}\}_{1\leq
j\leq m}$ given by Corollary 6.1 shows that

\bigskip

\noindent \textbf{Corollary 6.2. }\textsl{For the five exceptional Lie
groups one has}

\begin{quote}
$\mathcal{A}(G_{2})=\mathbb{Z}[y_{3}]/\left\langle
2y_{3},y_{3}^{2}\right\rangle $;

$\mathcal{A}(F_{4})=\mathbb{Z}[y_{3},y_{4}]/\left\langle
2y_{3},y_{3}^{2},3y_{4},y_{4}^{3}\right\rangle $;

$\mathcal{A}(E_{6})=\mathbb{Z}[y_{3},y_{4}]/\left\langle
2y_{3},y_{3}^{2},3y_{4},y_{4}^{3}\right\rangle $;

$\mathcal{A}(E_{7})=\mathbb{Z}[y_{3},y_{4},y_{5},y_{9}]/\left\langle
2y_{3},3y_{4},2y_{5},2y_{9},y_{3}^{2},y_{4}^{3},y_{5}^{2},y_{9}^{2}\right%
\rangle $;

$\mathcal{A}(E_{8})=\mathbb{Z}[y_{3},y_{4},y_{5},y_{6},y_{9},y_{10},y_{15}]/$

$\quad \left\langle
2y_{3},3y_{4},2y_{5},5y_{6},2y_{9},3y_{10},2y_{15},y_{3}^{8},y_{4}^{3},y_{5}^{4},y_{9}^{2},2y_{6}^{5}-y_{10}^{3}+y_{15}^{2}\right\rangle 
$.$\square $
\end{quote}

\noindent \textbf{Proof of Theorem 1.3. }Let $G$ be a simple Lie group with
rank $n$. By Theorem 1.2 the numbers of generators and relations in the
presentation of the ring $H^{\ast }(G/T)$ in (1.2) for $G\neq E_{8}$, and in
(1.3) for $G=E_{8}$, are both $n+m$. We shall show that, \textsl{without the
constraint that the generating set }$\{x_{1},\cdots ,x_{k}\}$\textsl{\
consists of Schubert classes, }this number is minimum with respect to any
presentation of the ring $H^{\ast }(G/T)$ in the form (1.1).

Let $\{s_{1},\cdots ,s_{h}\}\subset H^{\ast }(G/T)$ be a subset that
generates the ring $H^{\ast }(G/T)$ multiplicatively. Since the set $\left\{
\omega _{1},\cdots ,\omega _{n}\right\} $ of fundamental weights is a basis
of the group $H^{2}(G/T)$ and none of which can be expressed as a polynomial
in the lower degree ones, we can assume that $h\geq n$ and that $%
s_{i}=\omega _{i}$ for $1\leq i\leq n$. Now (6.4) implies that the quotient
ring $\mathcal{A}(G)$ is generated by those $p(s_{j})$ with $n+1\leq i\leq h$%
. Since $m$ is the minimal number of generators required to present the
quotient ring $\mathcal{A}(G)$ by (6.5), we have further that $h-n\geq m$.
This shows that $m+n$ is the least number of a set of generators of the ring 
$H^{\ast }(G/T)$. In particular $h(G,T)=n+m+1$.

To show $n+m$ is the least number of relations to characterize $H^{\ast
}(G/T)$ we can assume, by the remark after Definition 1.1, that $\left\{
h_{1},\cdots ,h_{q}\right\} $ is a set of homogeneous polynomials in $%
\left\{ \omega _{i},y_{d_{j}}\right\} _{1\leq i\leq n,1\leq j\leq m}$ which
satisfies that

\begin{enumerate}
\item[(6.6)] $H^{\ast }(G/T)=\mathbb{Z}\left[ \omega _{i},y_{d_{j}}\right]
_{1\leq i\leq n,1\leq j\leq m}/\left\langle h_{1},\cdots ,h_{q}\right\rangle 
$.
\end{enumerate}

\noindent Then one gets in addition to (6.5) another presentation of the
quotient

\begin{quote}
$\mathcal{A}(G)=\mathbb{Z}\left[ y_{d_{1}},\cdots ,y_{d_{m}}\right] _{1\leq
j\leq m}/\left\langle \overline{h}_{1},\cdots ,\overline{h}_{q}\right\rangle 
$, $\overline{h}_{i}=h_{i}\mid _{\omega _{1}=\cdots =\omega _{n}=0}$.
\end{quote}

\noindent Comparing this with (6.5) and in view of the sets $%
\{d_{j}\}_{1\leq j\leq m}$, $\{k_{j}\}_{1\leq j\leq m}$ of integers given by
Corollary 6.1, we can assume further that $q\geq m$ and that $\overline{h}%
_{i}=p_{i}y_{d_{i}}$ for all $1\leq i\leq m$, where the latter is equivalent
to

\begin{enumerate}
\item[(6.7)] $h_{i}=p_{i}y_{d_{i}}+\gamma _{i}$ with $\gamma _{i}\in
\left\langle \omega _{1},\cdots ,\omega _{n}\right\rangle ,$ $1\leq i\leq m$.
\end{enumerate}

\noindent Since $H^{\ast }(G/T;\mathbb{Q})=H^{\ast }(G/T)\otimes \mathbb{Q}$
one gets by (6.6) and (6.7) that

\begin{enumerate}
\item[(6.8)] $H^{\ast }(G/T;\mathbb{Q})=\mathbb{Q}[\omega _{1},\cdots
,\omega _{n}]/\left\langle \widetilde{h}_{m+1},\cdots ,\widetilde{h}%
_{q}\right\rangle $,
\end{enumerate}

\noindent where $\widetilde{h}_{t}$, $m+1\leq t\leq q$, is the polynomial
obtained from $h_{t}$ by substituting in $y_{d_{i}}=-\frac{1}{p_{i}}\gamma
_{i}$ by the relations (6.7). Since the variety $G/T$ is of finite
dimensional we must have $\dim H^{\ast }(G/T;\mathbb{Q})<\infty $.
Consequently $q-m\geq n$ by (6.8). That is, in the presentation (6.6) one
must have $q\geq n+m$. This completes the proof.$\square $

\bigskip

Let $D(\mathcal{A}(G))$ be the ideal of decomposable elements of the ring $%
\mathcal{A}(G)$, and let $q:$ $\mathcal{A}(G)\rightarrow \overline{\mathcal{A%
}}(G):=\mathcal{A}(G)/D(\mathcal{A}(G))$ be the quotient map. In view of
(6.5) the graded group $\overline{\mathcal{A}}(G)$ is determined by the data 
$\{d_{j};p_{j}\}_{1\leq j\leq m}$ as

\begin{enumerate}
\item[(6.9)] $\overline{\mathcal{A}}(G)=\mathbb{Z}\dbigoplus\limits_{1\leq
j\leq m}\overline{\mathcal{A}}^{d_{j}}(G)$ with $\overline{\mathcal{A}}%
^{d_{j}}(G)=\mathbb{Z}_{p_{j}}$.
\end{enumerate}

\noindent The proof of Theorem 1.3 is applicable to show that

\bigskip

\noindent \textbf{Theorem 6.3. }\textsl{A set }$S=\left\{ x_{d_{1}},\cdots
,x_{d_{m}}\right\} $\textsl{\ of Schubert classes on }$G/T$ \textsl{is
special if and only if the class }$q\circ p(x_{d_{j}})$ \textsl{is a
generator of the cyclic group} $\overline{\mathcal{A}}^{2d_{j}}(G)$,\textsl{%
\ }$1\leq j\leq m$.$\square $

\bigskip

We give an application of Theorem 1.2. For a Lie group $G$ with a maximal
torus $T$ consider corresponding fibration $\pi :G\rightarrow G/T$. In \cite%
{G} Grothendieck introduced the \textsl{Chow ring} $A(G^{c})$ for the
reductive algebraic group $G^{c}$ corresponding to $G$, and proved the
relation

\begin{enumerate}
\item[(6.10)] $A(G^{c})=$ Im $\{\pi ^{\ast }:H^{\ast }(G/T)\rightarrow
H^{\ast }(G)\}$.
\end{enumerate}

\noindent On the other hand resorting to the Leray--Serre spectral sequence
of $\pi $ one can show that

\begin{enumerate}
\item[(6.11)] Im$\pi ^{\ast }=H^{\ast }(G/T)/\left\langle \func{Im}\tau
\right\rangle $ (see \cite[Lemma 4.3]{D0}),
\end{enumerate}

\noindent where $\tau :H^{1}(T)\rightarrow H^{2}(G/T)$ be the transgression
in the fibration $\pi $ \cite[p.185]{Mc}. Granted with the explicit
presentation of the rings $H^{\ast }(G/T)$, as well as the formula \cite[%
formula (3.4)]{D0} for $\tau $, formula (6.11) is ready to apply to yield
formulae for the ring $A(G^{c})$ by Schubert classes on $G/T$.

As examples, if $G$ is $1$--connected, then $\left\langle \func{Im}\tau
\right\rangle =\left\langle \omega _{1},\cdots ,\omega _{n}\right\rangle $
by \cite[formula (3.4)]{D0}. Formula (6.11) implies that

\begin{quote}
$A(G^{c})=\mathcal{A}(G)$ (see Corollary 6.2)
\end{quote}

\noindent Similarly, for the adjoint Lie groups $PG$ with $G=SU(n)$, $%
Sp(n),E_{6}$ and $E_{7}$ one has (see \cite[formula (6.2)]{D0})

\begin{enumerate}
\item[(6.12)] 
\begin{tabular}{l}
$A(PSU(n)^{c})=\frac{\mathbb{Z}[\omega _{1}]}{\left\langle b_{r}\omega
_{1}^{r}\mid 1\leq r\leq n\right\rangle }$ with $b_{r}=g.c.d\{C_{n}^{1},%
\cdots ,C_{n}^{r}\}$, \\ 
$A(PSp(n)^{c})=\frac{\mathbb{Z}[\omega _{1}]}{\left\langle 2\omega
_{1},\omega _{1}^{2^{r+1}}\right\rangle }$, $n=2^{r}(2s+1)$, \\ 
$A(PE_{6}^{c})=\frac{\mathbb{Z}[\omega _{1},y_{3}^{\prime },y_{4}]}{%
\left\langle 3\omega _{1},2y_{3}^{\prime },3y_{4},x_{6}^{\prime 2},\omega
_{1}^{9},y_{4}^{3}\right\rangle }$, $y_{3}^{\prime }=y_{3}+\omega _{1}^{3}$,
\\ 
$A(PE_{7}^{c})=\frac{\mathbb{Z}[\omega _{2},y_{3},y_{4},y_{5},y_{9}]}{%
\left\langle 2\omega _{2},\omega _{2}^{2},2y_{3},3y_{4},2y_{5},2{y_{9},}%
y_{3}^{2},y_{4}^{3},y_{5}^{2},y_{9}^{2}\right\rangle }$.%
\end{tabular}
\end{enumerate}

\noindent \textbf{Remark 6.4.} For $G=Spin(n),G_{2}$ and $F_{4}$ Marlin \cite%
{M} obtained the ring $A(G^{c})$ by Schubert classes on $G/T$. For the
simple Lie groups Ka\v{c} \cite{K} computed the algebras $A(G^{c})\otimes 
\mathbb{F}_{p}$ with generators specified by the degrees.$\square $

\bigskip

\noindent \textbf{Remark 6.5.} For the earlier works studying the
presentation of the ring $H^{\ast }(G/T)$, see \cite{B,BS,M,T,TW,N1,N2}. We
remark that a basic requirement of intersection theory \cite{F} is to
present the cohomology $H^{\ast }(X)$ of a projective variety $X$ by
explicit described geometrical cycles, such as the Schubert classes of flag
manifolds, so that the intersection multiplicities can be computed by the
cup products on the ring $H^{\ast }(X)$. In this regard the approaches due
to Bott--Samelson \cite{BS} and Marlin \cite{M} are inspiring.

\bigskip

The authors would like to thank their referees for valuable suggestions. In
particular, the question of giving formulae for the Chow rings for the
adjoint Lie groups $PG$ has motivated the subsequent work \cite{D0} on the
integral cohomology of all compact Lie groups.


\begin{thebibliography}{99}
\bibitem{AH} M. Atiyah, F. Hirzebruch, Vector bundles and homogeneous
spaces, 1961 Proc. Sympos. Pure Math., Vol. III pp. 7--38.

\bibitem{BGG} I. N. Bernstein, I. M. Gel'fand, S. I. Gel'fand, Schubert
cells and cohomology of the spaces G/P, Russian Math. Surveys 28 (1973),
1-26.

\bibitem{B} A. Borel, Sur la cohomologie des espaces fibres principaux et
des espaces homogenes de groupes de Lie compacts, Ann. math. 57(1953),
115--207.

\bibitem{BH} A. Borel, F Hirzebruch, Characteristic classes and homogeneous
space I, Amer. J. Math., vol 80 (1958), 458-538.

\bibitem{BS} R. Bott, H. Samelson, The cohomology ring of $G/T$, Nat. Acad.
Sci. 41(1955), 490-492.

\bibitem{Ch} C. Chevalley, Sur les D\'{e}compositions Cellulaires des
Espaces G/B, in Algebraic groups and their generalizations: Classical
methods, W. Haboush ed. Proc. Symp. in Pure Math. 56 (part 1) (1994), 1-26.

\bibitem{CP} P. E. Chaput, and N. Perrin, Towards a Littlewood-Richardson
rule for Kac-Moody homogeneous spaces, J. Lie Theory 22 (2012), no. 1,
17--80.

\bibitem{De} M. Demazure, D\'{e}singularisation des vari\'{e}t\'{e}s de
Schubert g\'{e}n\'{e}ralis\'{e}es, Ann. Sci. \'{E}cole Norm. Sup. (4) 7
(1974), 53--88.

\bibitem{D} H. Duan, Multiplicative rule of Schubert classes, Invent. Math.
159(2005), no. 2, 407--436; 177(2009), no.3, 683--684.

\bibitem{D0} H. Duan, The cohomology of compact Lie groups, math.AT.
arXiv:1502.00410.

\bibitem{D1} H. Duan, Self-maps of the Grassmannian of complex structures.
Compositio Math. 132 (2002), no. 2, 159--175.

\bibitem{DL} H. Duan, SL. Liu, The isomorphism type of the centralizer of an
element in a Lie group, Journal of algebra, 376(2013), 25-45.

\bibitem{DZ0} H. Duan and Xuan Zhao, The classification of cohomology
endomorphisms of certain flag manifolds, Pacific J. Math. 192 (2000), no. 1,
93--102.

\bibitem{DZ1} H. Duan, Xuezhi Zhao, Erratum: Multiplicative rule of Schubert
classes, Invent. Math. 177 (2009), 683--684.

\bibitem{DZ2} H. Duan, Xuezhi Zhao, The Chow rings of generalized
Grassmannians, Found. Math. Comput. Vol.10, no.3(2010), 245--274.

\bibitem{DZ3} H. Duan, Xuezhi Zhao, Schubert calculus and cohomology of Lie
groups, math.AT (math.AG). arXiv:0711.2541.

\bibitem{DZ4} H. Duan, Xuezhi Zhao, Schubert calculus and the Hopf algebra
structures of exceptional Lie groups, Forum. Math. Vol.26, no.1(2014),
113-140.

\bibitem{DZ5} H. Duan and Xuezhi Zhao, Algorithm for multiplying Schubert
classes, International Journal of Algebra and Computation, vol 16,
No.6(2006), 1197-1210.

\bibitem{F} W. Fulton, Intersection theory, Springer--Verlag, 1998.

\bibitem{G} A. Grothendieck, Torsion homologique et sections rationnelles,
Sem. C. Chevalley, ENS 1958, expos\'{e} 5, Secreatariat Math. IHP, Paris,
1958.

\bibitem{H} J. E. Humphreys, Introduction to Lie algebras and representation
theory, Graduated Texts in Math. 9, Springer-Verlag New York, 1972.

\bibitem{Hus} D. Husemoller, Fibre bundles, Second edition. Graduate Texts
in Mathematics, No. 20. Springer-Verlag, New York-Heidelberg, 1975.

\bibitem{K} V.G. Ka\v{c}, Torsion in cohomology of compact Lie groups and
Chow rings of reductive algebraic groups, Invent. Math. 80(1985), no. 1,
69--79.

\bibitem{Kl} S. Kleiman, Intersection theory and enumerative geometry: A
decade in review, Algebraic geometry, Proc. Summer Res. Inst.,
Brunswick/Maine 1985, part 2, Proc. Symp. Pure Math. 46, 321-370 (1987).

\bibitem{Ku} S. Kumar, Kac-Moody groups, their flag varieties and
representation theory, Progress in Mathematics, 204. Birkha\"{u}ser Boston,
Inc., Boston, MA, 2002.

\bibitem{LG} V. Lakshmibai, N. Gonciulea, The Flag variety,
Hermann-Actualites Mathematiques (2001).

\bibitem{M} R. Marlin, Anneaux de Chow des groupes alg\'{e}riques $SU(n)$, $%
Sp(n)$, $SO(n)$, $Spin(n),G_{2},F_{4}$, C. R. Acad. Sci. Paris, A 279(1974),
119--122.

\bibitem{Mc} J. McCleary, A user's guide to spectral sequences, Second
edition. Cambridge Studies in Advanced Mathematics, 58. Cambridge University
Press, Cambridge, 2001.

\bibitem{N1} M. Nakagawa, The integral cohomology ring of $E_{7}/T$, J.
Math. Kyoto Univ. 41 (2001), 303-321.

\bibitem{N2} M. Nakagawa, The integral cohomology ring of $E_{8}/T$, Proc.
Japan Acad. Ser. A Math. Sci. 86 (2010), 64-68.

\bibitem{Sch} H. Schubert, Kalk\"{u}l der abz\"{a}hlenden Geometrie, Berlin,
Heidelberg, New York: Springer-Verlag (1979).

\bibitem{So} F. Sottile, Four entries for Kluwer encyclopaedia of
Mathematics, arXiv: Math. AG/0102047.

\bibitem{T} H. Toda, On the cohomology ring of some homogeneous spaces, J.
Math. Kyoto Univ. 15(1975), 185--199.

\bibitem{TW} H. Toda, T. Watanabe, The integral cohomology ring of $F_{4}/T$
and $E_{6}/T$. J. Math. Kyoto Univ. 14 (1974), 257-286

\bibitem{Wa0} B. L. van der Waerden, Topologische Begr\"{u}ndung des Kalk%
\"{u}ls der abz\"{a}hlenden Geometrie, Math. Ann. 102 (1930), no. 1,
337--362.

\bibitem{Wa} B. L. van der Waerden, The Foundation of Algebraic Geometry
from Severi to Andr\'{e} Weil, Archive for History of Exact Sciences, Vol.
7, No. 3 (26.V.1971), pp. 171-180.

\bibitem{W} A. Weil, Foundations of algebraic geometry, American
Mathematical Society, Providence, R.I. 1962.
\end{thebibliography}
\end{document}